\documentclass{amsart}\usepackage{amssymb,euscript,amsfonts,bbm,graphicx}%\usepackage[mathscr]{},
\usepackage{enumerate}
\newtheorem{Lemma}{Lemma}[section]\newcommand{\bel}{\begin{Lemma}}\newcommand{\eel}{\end{Lemma}}
\newtheorem{Example}[Lemma]{Example}\newcommand{\bex}{\begin{Example}\rm}\newcommand{\eex}{\end{Example}}
\newtheorem{Proposition}[Lemma]{Proposition}\newcommand{\bprop}{\begin{Proposition}}\newcommand{\eprop}{\end{Proposition}}
\newtheorem{Definition-Proposition}[Lemma]{Definition-Proposition}

\def\bpr{~\\{\em Proof.\ }}\newcommand{\epr}{$\bull$\\}
\newtheorem{Theorem}[Lemma]{Theorem}\newcommand{\bthe}{\begin{Theorem}}\newcommand{\ethe}{\end{Theorem}}
\newtheorem{Definition}[Lemma]{Definition}\newcommand{\bed}{\begin{Definition}}\newcommand{\eed}{\end{Definition}}
\newtheorem{Remark}[Lemma]{Remark}\newcommand{\beR}{\begin{Remark}\rm}\newcommand{\eeR}{\end{Remark}}
\newtheorem{Corollary}[Lemma]{Corollary}\newcommand{\bcor}{\begin{Corollary}}\newcommand{\ecor}{\end{Corollary}}

\newcommand{\beq}{\begin{equation}}\newcommand{\eeq}{\end{equation}}
\newcommand{\beqn}{\begin{equation*}}\newcommand{\eeqn}{\end{equation*}}
\newcommand{\bem}{\begin{displaymath}}\newcommand{\eem}{\end{displaymath}}
\newcommand{\beqa}{\begin{eqnarray}}\newcommand{\eeqa}{\end{eqnarray}}
\newcommand{\bee}{\begin{enumerate}}\newcommand{\eee}{\end{enumerate}}
\newcommand{\bei}{\begin{itemize}}\newcommand{\eei}{\end{itemize}}
\newcommand{\bet}{\begin{tabular}{cccccccc}}\newcommand{\eet}{\end{tabular}}
\newcommand{\bpm}{\begin{pmatrix}}\newcommand{\epm}{\end{pmatrix}}
\newcommand{\bM}{\begin{matrix}}\newcommand{\eM}{\end{matrix}}
\newcommand{\ber}{\begin{array}{l}}\newcommand{\eer}{\end{array}}

\def\o{o}\newcommand{\frD}{\mathfrak{D}}

\def\bull{\vrule height .9ex width .9ex depth -.1ex }
\newcommand{\quotient}[2]{{\left.\raisebox{1.6ex}{$#1$}\!\!\!\!\!{\scalebox{2}{\ensuremath\diagup}}
\!\!\!\!\!\raisebox{-1ex}{$#2$}\right.}}

\newcommand{\quotients}[2]{{\footnotesize\left.\raisebox{0.4ex}{$#1$}\! / \!\raisebox{-0.4ex}{$#2$}\right.}}

\def\di{\partial}
\def\bl{\langle}\def\br{\rangle}\def\ra{\rightarrow}
\def\isom{\xrightarrow{\sim}}
\def\into{\stackrel{i}{\hookrightarrow}}\def\proj{\stackrel{\pi}{\ra}}

\def\cN{\mathcal{N}}\def\cO{\mathcal{O}}\def\cp{{\mathfrak{p}}}

\def\cU{\mathcal{U}}
\def\cX{\mathcal{X}}\def\cZ{\mathcal{Z}}
\def\cm{{\frak m}}

\def\C{\mathbb{C}}%\def\C{\C}

\def\k{\mathbbm{k}}\def\K{\mathbb{K}}\def\P{\mathbb{P}}
\def\R{\mathbb{R}} %\def\R{\mathbb{R}}

\def\Z{\mathbb{Z}}

\def\al{\alpha}\def\be{\beta}\def\De{\Delta}\def\ga{\gamma}\def\Ga{\Gamma}
\def\ep{\epsilon}\def\la{\lambda}\def\La{\Lambda}\def\si{\sigma}
\def\Si{\Sigma}\def\Om{\Omega}

\def\tf{{\tilde{f}}}
\def\tI{{\tilde{I}}}

  \def\tX{{\tilde{X}}}

\def\ua{\underline{a}}
\def\uf{\underline{f}}\def\ug{\underline{g}}\def\um{{\underline{m}}}\def\us{\underline{s}}
\def\ux{\underline{x}}\def\uy{\underline{y}}\def\uz{\underline{z}}

\def\hi{{\hat{i}}}\def\hj{{\hat{j}}}

\def\empty{\varnothing}
\def\suml{\sum\limits}\def\oplusl{\mathop\oplus\limits}
\def\liml{\lim\limits}
\def\capl{\mathop\cap\limits}

\def\prodl{\prod\limits}

\def\lSi{{\overline{\Sigma}}}

\def\uI{{\underline{I}}}

\def\dv{{\vec{\di}}}

\def\smin{\setminus}\def\sset{\subset}\def\sseteq{\subseteq}

\def\omp{ordinary multiple point}\def\Db{{\De^\bot}}

\def\?{{\bf ???}}\def\cOn{{\cO_{(\k^N,\o)}}}
\newcommand{\bin}[2]{\binom{#1}{#2}}

\title{D\MakeLowercase{iscriminant of the ordinary transversal singularity type.} T\MakeLowercase{he local aspects.}}
\author{D\MakeLowercase{mitry} K\MakeLowercase{erner} \MakeLowercase{and} A\MakeLowercase{ndr\'as} N\MakeLowercase{\'emethi}}
\date{\today}

\address{Department of Mathematics, Ben Gurion University of the Negev, Israel}
\email{dmitry.kerner@gmail.com}

\address{Alfr\'ed R\'enyi Institute of Mathematics,
Hungarian Academy of Sciences,
Re\'altanoda utca 13-15, H-1053, Budapest, Hungary \newline
 \hspace*{4mm} ELTE - University of Budapest, Dept. of Geometry, Budapest, Hungary \newline \hspace*{4mm}
BCAM - Basque Center for Applied Math.,
Mazarredo, 14 E48009 Bilbao, Basque Country – Spain}
\email{nemethi.andras@renyi.mta.hu }
\thanks{D.K. was supported by the grant FP7-People-MCA-CIG, 334347.
Part of the work was done during D.K.'s posdoctoral fellowship in the University of Toronto.}
\thanks{A.N. was partially supported by NKFIH Grant  112735 and
ERC Adv. Grant LDTBud of A. Stipsicz at R\'enyi Institute of Math., Budapest}
\thanks{We thank V. Goryunov, P. Milman, D. Siersma,  D. van Straten, B. Sturmfels for important advices.}

\begin{document}\setcounter{secnumdepth}{6} \setcounter{tocdepth}{1}

\begin{abstract}
Consider a space $X$ with the singular locus, $Z=Sing(X)$, of positive dimension. Suppose both $Z$ and $X$ are locally complete intersections.
 The  transversal type of $X$ along $Z$ is generically constant but at some points of $Z$ it degenerates. We introduce (under certain conditions)
 the discriminant of the transversal type, $\Db$, a subscheme of $Z$, that reflects these degenerations
whenever the generic transversal type is `ordinary'.

The scheme structure of $\Db$ is imposed by various compatibility properties and is often non-reduced. We establish the basic properties of $\Db$:
 it is a Cartier divisor in $Z$, functorial under base change,  flat under some deformations of $(X,Z)$, and compatible with pullback under some morphisms, etc.
  Furthermore, we study the local  geometry of $\Db$, e.g. we compute its multiplicity at a point, and we obtain the resolution
  of $\cO_\Db$ (as $\cO_Z$-module) and
  study the locally defining equation.
\end{abstract}
\maketitle\tableofcontents

\section{Introduction}

\subsection{The setup}\label{Sec.Setup}
Let $\k$ be an algebraically closed field of zero characteristic, e.g. $\k=\C$. Let $M$ be either a
smooth irreducible algebraic variety (over $\k$), or, for $\k=\C$, a complex-analytic connected manifold.
Let $X\sset M$ be a reduced subscheme with non-isolated singularities. We assume that $Z:=Sing(X)$ is connected,
 otherwise one fixes a connected component $Z\sset Sing(X)$ and replaces $X$ by some neighborhood of $Z$.
 We always take $Z$ with its reduced  structure.

\

In many examples of non-isolated singularities one observes the following pattern.
For each smooth point $\o\in Z$ consider a smooth germ, $(L^\bot,\o)\sset(M,\o)$,  transversal to $(Z,\o)$, such that $(L^\bot,\o)\cap(Z,\o)=\{\o\}$.
The singularity $(L^\bot\cap X,\o)$ is usually isolated and its type
 is in some sense generically constant along $Z$ (thus it is called the ``transversal singularity type").
 The points where the transversal singularity type degenerates usually form a subset of codimension 1 in $Z$.
  It is natural to call this subset the {\em discriminant of the transversal type}, $\Db\sset Z$. This is the target of our work.

With the following examples we try to give some intuition and the guiding principles. The precise discussion will be given later.

First we show that at {\em some points} the transversal type is not well defined.
\bex\label{Ex.Transversal.Type.Depends.choice.of.Section}
Consider the singular surface $X=\{x^2z=y^2\}\sset\k^3$. Its (reduced) singular locus is  the
line $Z=\{x=y=0\}\sset\k^3$. This is the classical Whitney umbrella/pinch point/$D_\infty$ point.
 For the generic point $\o\in Sing(X)$,
i.e. for $z\neq0$, the transversal singularity, $(X,\o)\cap(L^\bot,\o)$, is the plane curve singularity of type $A_1$, i.e. two smooth non-tangent branches.
 As $z\to0$ the transversal singularity degenerates, at the origin the transversal type is not well defined.
 Indeed,  we choose the transversal section $(L^\bot,\o)$ among those defined by equation  $z=ax+by+(higher\ order\ terms)$.
\bei
\item For $a\neq0$ the intersection $(X,\o)\cap(L^\bot,\o)$ is a cusp, $A_2$.
\item For $(L^\bot,\o)=\{z=x^{n-1}\}$ the intersection $(X,\o)\cap(L^\bot,\o)$ is $A_n$-singularity.
\item For $(L^\bot,\o)=\{z=0\}$ the intersection $(X,\o)\cap(L^\bot,\o)$ is a double line,  a non-isolated singularity.
\eei

Therefore the expectation is that the point $(0,0,0)$ belongs to the discriminant $\Db$.
\eex

\

  The following example suggests that sometimes the scheme structure on $\Db$ should be taken non-reduced.
\bex\label{Ex.Whitney.Umbrella}
Let $X=\{x^2z^q=y^2+x^3\}\sset\k^3$ for $q\ge1$. As before, the singular locus is $Z=\{x=y=0\}$ and the transversal type degenerates as $z\neq0$.
Consider the deformation: $X_t=\{x^2(z^q-t)=y^2+x^3\}\sset\k^3$ for  $t\in (\k,\o)$. It preserves the singular
 set: $Sing(X_t)=\{x=0=y\}$. For $t\neq0$ the discriminantal point $(0,0,0)$ splits into $q$ points $\{x=y=z^q-t=0\}$, each of them being of $D_\infty$-type.
 Thus, for $t=0$, it is natural
to consider  the point $(0,0,0)\in\Db$ with multiplicity $q$ (or a multiple of $q$). One can say
roughly that for $q>1$ the transversal type degenerates (as $z\to0$) `faster'.
(In examples of \S\ref{Sec.Db.Definition} we give
 other reasons for non-reducedness of $\Db$.)
\eex
We remark that the naive geometric consideration of the transversal section, $(X,\o)\cap(L^\bot,\o)$, does not work at the singular points of $Z$.
 Hence it should be replaced by an algebraic counterpart.%, see \S\ref{Sec.Background.Normal.Cone} and \S\ref{Sec.Background.Ordinary.Singularities}.

\subsection{Assumptions}\label{Sec.Intro.Assumptions}
The definition of transversal type and its discriminant in the full generality seems out of reach at the present stage.
Indeed, this would use the equisingularity theory in arbitrary dimension and codimension for arbitrary classes of singularities
 (e.g. whenever $(Z,\o)$ is not necessarily Gorenstein or Cohen-Macaulay). Thus we work under the following assumptions.
  (The precise definition, examples and properties are in \S\ref{Sec.Background}.)
\bei
\item The (reduced) singular locus, $Z=Sing(X)$ or $(Z,\o)=Sing(X,\o)$, is a locally complete intersection at each point (l.c.i.).
\item For each point $\o\in Z$ the germ $(X,\o)$ is  a {\em strictly complete intersection} over $(Z,\o)$ (s.c.i.).
This is a strengthening of the notion of complete intersection, needed to ensure that the strict transform under blowup along $Z$ is again
a complete intersection.  %, see \S\ref{Sec.Background.Strictly.Complete.Intersections}.
 In particular, if $Z$  has several irreducible components then the   multidegree  of $\P T_{(L^\bot\cap X,\o)}$ at
{\em generic} points of each  component is the same.
\item The transversal type of $X$ along $Z$  is  generically `ordinary'. Namely,
 for sufficiently generic point $\o\in Z$, the projectivized tangent cone, $\P T_{(L^\bot\cap X,\o)}$,
 is a smooth complete intersection of expected dimension.
\eei

\

Under these assumptions we define the discriminant of transversal type, $\Db=\Db_{X/Z}$, (with the natural scheme structure), and establish
 some local and global properties. (The further global properties are established in
 \cite{Kazarian.Kerner.Nemethi}.)
 For the history of the question and some known results see \S\ref{Sec.Intro.History}.

\subsection{On the choice of scheme structure of the discriminant}\label{Sec.Intro.Scheme.Structure}
In simple cases, like that of examples \ref{Ex.Transversal.Type.Depends.choice.of.Section}, \ref{Ex.Whitney.Umbrella}, it
  is obvious which points belong to $\Db$. This determines $\Db\sset Z$ as a subset,
not as a subscheme. On the other hand, it is less obvious whether/when the singular points  of $Z$ belong to $\Db$.

Our definition of the subscheme $(\Db,\o)\sset(Z,\o)$ is guided by the wish-list of the following natural properties:
\bee
\item ({\em Normalization}) For the classical Whitney umbrella, $\{x^2z=y^2\}\sset\k^3$, the discriminant is the reduced point $(0,0,0)\in\k^3$.
More generally, for a $D_\infty$ point, $\{x_0x^2_1+\suml^n_{i=2}x^2_i=0\}\sset(\k^{n+1},0)$, the
  discriminant is the reduced point $(0,\dots,0)\in (Z,\o)=\{x_1=\cdots=x_n=0\}\sset\k^{n+1}$.
 Even more generally, suppose the germ $(Z,\o)$ is smooth and the multiplicity of $X$ along $Z$ is locally constant
  at $\o$. Take the generic section $(L^\bot\cap X,\o)$, suppose the projectivization of the tangent cone,
  $\P T_{(L^\bot\cap X,\o)}$, has just one $A_1$ singularity. Then $\Db\sset Z$ is reduced at $\o$.
\item ({\em Behaviour in families.} Note that the family $\{\Db(X_t)\}$ of example \ref{Ex.Whitney.Umbrella} is flat.)
  Suppose that a flat family $\cX=\{X_t\}_{t\in(\k^1,\o)}\to(\k^1,\o)$ satisfies:
\bei
\item
the family $\cZ=\{Z_t=Sing(X_t)\}_{t\in(\k^1,0)}\to(\k^1,\o)$ is flat;
\item the generic multiplicity of $X_t$ along $Z_t$ does not vary with $t$;
\item for any $t$ the transversal type of $X_t$ along $Z_t$ is generically ordinary (see \S\ref{Sec.Intro.Assumptions}).
\eei
Then the family $\{\Db(X_t)\}_{t\in(\k^1,0)}$ is flat.

\item ({\em Pullback of the classical discriminant}) Suppose $Z$ (or its germ at a point) is smooth. Take the strict transform under blowup, $Bl_Z(M)\supset\tX$.
 The exceptional divisor, $E\sset Bl_Z(M)$ induces the family of projective complete intersections, $\tX\cap E\to Z$. Thus one has a (rational)
  map from $Z$ to the parameter space of projective complete intersections. (In the hypersurface case this parameter space is $|\cO_{\P^n}(d)|$,
   for complete intersections one can take e.g. $\prod|\cO_{\P^n}(d_i)|$.) In this parameter space we have the classical discriminant $\De$. Then $\Db$
    should be the pullback of $\De$.

\item
({\em The image of the critical locus}) Suppose the fibres of the projection $\tX\cap E\stackrel{\pi}{\to}Z$ are (generically) of dimension $d$.
The critical locus, $Crit(\pi)\sseteq \tX\cap E$ is defined via the relative cotangent sheaf, $\Om^1_{\quotients{\tX\cap E}{Z}}$,
by the Fitting ideal $Fitt_d(\Om^1_{\quotients{\tX\cap E}{Z}})\sseteq\cO_{\tX\cap E}$.
 Then $\Db$ is the image of $Crit(\pi)$, with the Fitting scheme structure, $I_{\Db/Z}=Fitt_0(\pi_*\cO_{Crit\pi})$.
\eee

\

We define the subscheme $\Db\sset Z$ in \S\ref{Sec.Db.Definition}, it has all these properties.

\subsection{Additional basic properties of $\Db$}
Besides the minimal requirements listed above, the discriminant of transversal singularity type possesses (as a scheme)
 various other  nice/natural properties.

\bee
\item The scheme structure of $\Db$ is completely determined
 by the `infinitesimal neighborhood' of $Sing(X)$ in $X$, more precisely, by the exceptional
divisor of blowup:
$(E,\tX\cap E)\sset (Bl_ZM,\tX)$, see \S\ref{Sec.Db.Gen.Prop.Defined.by.Infinit.Neighb}. In this way it is independent of those `higher-order'
 degenerations of
$(L^\bot\cap X,\o)$ that preserve the tangent cone.
(In particular we do not see any direct relation of $\Db$ to the L\^{e} cycles of \cite{Massey-book}, see \S\ref{Sec.Db.vs.Le.Numbers}.)

\item
 (The discriminant pulls back.) Given a morphism $M_1\stackrel{\phi}{\to}M_2$, inducing $X_1=\phi^*(X_2)$ and $Z_1=\phi^*(Z_2)$.
  Suppose $Z_i$ are reduced l.c.i. and are connected components of $Sing(X_i)$. Suppose $X_i$ are s.c.i. over $Z_i$ at
   each point  and $X_i$ are generically ordinary along $Z_i$, with the same  multiplicity sequences. Then $\Db_{X_1/Z_1}=\phi^*\Db_{X_2/Z_2}$,
   see \S\ref{Sec.Db.Gen.Prop.Pulls.Back}.

A particular case of the statement is the following: given a smooth hypersurface germ $(M_1,\o)\sset(M_2,\o)$, such that the
 tangent cones intersection $T_{(M_1,\o)}\cap T_{(X,\o)}$ is generic enough, then $\Db_{M_1\cap X}=\Db_X \cap (M_1,\o)$.

\item
 The (local) defining equation of $\Db$ is obtained by elimination procedure and thus cannot be written explicitly in the full generality.
  Yet, following the tradition, we present the discriminant as the determinant of a matrix. More precisely, we establish the (traditional)
   free resolution of $\cO_{(\Db,\o)}$,    as a module over $\cO_{(Z,\o)}$, see \S\ref{Sec.Db.Resolution.Defining.Equation.Tangent.Cone}.
 We use this resolution to get some information about the monomials of the defining equation of $\Db$.
In particular,  in the weighted-homogeneous case, we compute the total (weighted) degree of monomials that occur in the discriminantal polynomial,
 see Proposition \ref{Thm.Total.Degree.Disriminant}.

\item Flatness of a deformation of $\Db$ (under the deformation of $X$, see (2) of \S\ref{Sec.Intro.Scheme.Structure}) means the following: the sheaves $\cO_{Z_t}(-\Db_t)$ glue to a locally free
 sheaf of ideals $I_{\{\Db_t\}}$ on $\cZ=\{Z_t\}$
 and the schemes $\Db_t$ glue to a Cartier divisor on $\cZ$.
 If $\{X_t\}$ are not equimultiple along $\{Z_t\}$ or the induced deformation $\{Z_t\}$ is not flat then the family $\Db_t$ is
  not flat and in general is not semi-continuous in any sense, see \S\ref{Sec.Db.Gen.Prop.Deforms.Flatly}.

\item (The multiplicity of $\Db$ at a point.) Given the projection  $\tX\cap E\proj Z$ suppose the fibre $\pi^{-1}(\o)$ has only isolated singularities.
 Then $(\Db,\o)=\sum  (\Db_i,\o)$, the sum of Cartier divisors corresponding to the singular points of $\pi^{-1}(\o)$. Thus it is enough to assume that
 $\pi^{-1}(\o)$ has only one singular point.  In the hypersurface case,
  suppose in some local coordinates $\uz$ on $(Z,\o)$ and $\ux$ on $\pi^{-1}(\o)$ the locally defining equation of $\pi^{-1}(\o)$ is $f(\ux)+g(\uz)$.
  Then $mult(\Db,\o)=\mu(\pi^{-1}(\o))\cdot mult(g(\uz))$. For complete intersection we obtain a similar result, using the L\^{e}-Greuel formula,
    see \S\ref{Sec.Db.Multiplicity}.
 \item In \S\ref{Sec.Db.Stratification.Singular.Locus} we define a further stratification of $\Db$, corresponding to the higher degenerations of transversal type.
\eee

\

We emphasize that in the hypersurface case most statements of our paper appear in the standard literature. But the case of complete intersections is less known.

\subsection{History and motivation}\label{Sec.Intro.History}
\bei
\item The discriminant of transversal singularity type appears naturally
in geometry and singularity theory and in some particular cases
was considered already by Salmon, Cayley, Noether and Zeuthen, see \cite{Piene1978}.
 One
context where it appears is the image of the generic map from a smooth $n$--fold
into $\P^{n+1}$. The image has non-isolated ordinary
singularities, \cite[page 111]{Mond-Pellikaan}, (not to be confused with the `ordinary transversal
type' used in this paper). The natural question is to understand
their degenerations, as one runs along the singular locus.

\item The class of $\Db$ for projective surface,
$X\sset\P^3$, with ordinary singularities goes back (probably) to
the early history. For a computation see  \cite{Piene1978} (among
various other invariants).

\item The case of one-dimensional singular locus, i.e. $Z$ is a curve, with the generic transversal type $A_1$,
 was thoroughly studied by Siersma, see e.g. \cite{Siersma2000}. The local degree of the discriminant, called also
 `the virtual number of $D_\infty$ points' was studied in \cite{Pellikaan-PhD}, \cite{Pellikaan} and \cite{de Jong}.
 In particular, the authors show pathological behavior when $Z$ is
 not a locally complete intersection. In \cite{de Jong-de Jong} the degree of $[\Db]$ is computed for the
 case $X\sset M$ is a projective hypersurface, $Z=Sing(X)$ is of (pure) dimension one and the generic transversal type is $A_1$.
For the review of various related result see \cite[\S I.4.6]{AGLV2}.
For the recent results and applications to real singularities see \cite{van Straten}.

We emphasize that in Pelikaan-de Jong's approach the scheme structure on the discriminant is compatible with flat deformations, \cite[\S2.5]{de Jong},
 and the discriminant is reduced for Whitney umbrella. These two conditions determine the scheme structure uniquely,
 therefore their and our scheme structures (for non-isolated singularities of surfaces) coincide.
In example \ref{Ex.Discriminant.Eq.for.generic.transversal.A1} we show this directly.

\item One often considers the singular locus with the scheme structure defined by Jacobian ideal, $Sing(X)^{(J_f)}$,
 \cite{Aluffi-1995}, \cite{Aluffi-2005}. The scheme $Sing(X)^{(J_f)}$ also reflects the degenerations of transversal type. We emphasize,
that this Jacobian scheme structure is incompatible with flat deformations and it differs from the scheme structure of our paper.
\eei

\section{Preliminaries}\label{Sec.Background}
For the general introduction to singularities see \cite{AGLV},  \cite{Dimca92}, \cite{Looijenga-book} and
\cite{Seade}.

\subsection{Local neighborhoods}
Working locally, we consider germs of spaces,
$(Z,\o)\sset(X,\o)\sset(\k^N,\o)$. These germs can be
 algebraic, analytic (for $\k=\C$), formal, etc., the category is specified by the (local) ring of regular functions.
 The ring $\cOn$ is a regular (Noetherian) local ring over a field of zero characteristic, e.g. one of the following:
 $\k[x_1,\dots,x_n]_{(\cm)}$ (localization of the affine ring), or $\k\{x_1,\dots,x_n\}$ (the ring of analytic power series, for $\k\sseteq\C$),
 or  $\k\bl x_1,\dots,x_n\br$ (the ring of algebraic power series), or $\k[[x_1,\dots,x_n]]$.
  For a subgerm $(X,\o)\sset(\k^N,\o)$ the local ring is the quotient by the defining ideal,   $\cO_{(X,\o)}=\quotients{\cO_{(\k^N,\o)}}{I_{(X,\o)}}$.
 In many cases the algebraic germs are `too large and rigid', e.g.
 when speaking of irreducible components or rectifying locally a smooth variety. In such cases
 we take henselization or completion (i.e. we pass to henselian or formal germs).

\

If the germ $(X,\o)$ is not algebraic/analytic then one cannot take its ``small enough representative", e.g. a formal germ has no closed
 points besides the base point. Yet, using the standard algebra-geometry dictionary the ideas/notions of ``working near the origin" are applicable.
 One just translates a geometric statement/condition into the algebraic one, e.g.:
\bei
\item ``the points of the subgerm $(Z,\o)\sset(X,\o)$ satisfy \dots" is replaced by ``the ideal $I_{(Z,\o)}\sset\cO_{(X,\o)}$ satisfies \dots"
\item ``generic points of $(Z,\o)\sset(X,\o)$ satisfy \dots" is replaced by ``the localization of $\cO_{(Z,\o)}$, $\cO_{(X,\o)}$,
  at the prime components of $I_{(Z,\o)}$ satisfies \dots"
\eei

\

We denote the maximal ideal in the local ring $R$ by $\cm_R$ (e.g. $\cm_{(X,\o)}$, $\cm_{(Z,\o)}$) or just by $\cm$.

\subsection{Multiplicity at a point, generic vanishing order and symbolic powers of ideals}
The (Taylor) order or multiplicity of an element $f$ in a local ring $(R,\cm)$ is defined as usual: $mult_R(f)=max\{k|\ f\in\cm^k\}$.
 More generally, the order of $f$ with respect to an ideal $J\sset R$ is  $ord_J(f)=max\{k|\ f\in J^k\}$.

 The multiplicity of a  germ $(X,\o)\sset(\k^N,\o)$ of pure dimension $n$ is defined as $dim_\k\quotients{\cO_{(X,\o)}}{(l_1,\dots,l_n)}$,
  where $(l_1,\dots,l_n)$ is the ideal generated by any $n$-tuple of generic elements of $\cm$.

\

Let the germ $(Z,\o)\sset(\k^N,\o)$ be reduced. An element $f\in \cO_{(\k^N,\o)}$ has {\em generic order $\ge m$ along $(Z,\o)$} if its
 (Taylor) order at smooth points of $(Z,\o)$ is $\ge m$.
 The general definition of this property goes via the notion of symbolic powers, as follows.
  (We replace $\cO_{(\k^N,\o)}$ by $R$ and $I_{(Z,\o)}$ by $J$.)

 Let $R$ be a Noetherian ring and $J\sset R$ a primary ideal, whose corresponding prime is $\cp$. The $m$'th symbolic power is defined as
 \beq
 J^{(m)}:=(J^m\cdot R_\cp)\cap R=\{f\in R|\ \exists s\not\in\cp:\ sf\in J^m\}.
 \eeq
If the ideal $J$ is not primary but radical, one takes the primary decomposition $J=\cap J_i$ and defines $J^{(m)}=\cap J^{(m)}_i$.
In the most general case the definition goes as follows (see definition 3.5.1 of \cite{Vasconcelos}).
 For any ideal in a Noetherian ring, $J\sset R$, take the decomposition: $J=J'\cap L$, where $J'$  is the intersection of
 the primary ideals associated with the minimal primes of $J$, while $L$ is the intersection of primary ideals corresponding to embedded primes of $J$.
  Then $J^{(n)}:=(J^n)'$.

\bed
We say that $f$   is generically of order $\ge m$ on all the components of $V(J)$ if $f\in J^{(m)}$.
\eed
For the explanation that $f\in J^{(m)}$ means this geometric condition see \cite[\S3.9]{Eisenbud-book}.

One has the obvious inclusion $I^{(m)}_{(Z,\o)}\supseteq I^m_{(Z,\o)}$ and this inclusion can be proper.
\bex
Let $(Z,\o)=\{xy=yz=xz=0\}\sset(\k^3,\o)$, thus $I_{(Z,\o)}=(x,y)\cap(y,z)\cap(x,z)$.
The generic order of $f=xyz$ along $(Z,\o)$ is $2$, thus $xyz\in I^{(2)}_{(Z,\o)}$, but $xyz\not\in I^2_{(Z,\o)}$.
 In fact we have the primary decomposition:
  $ I^2_{(Z,\o)}=(x,y)^2\cap (y,z)^2\cap (x,z)^2\cap (x^2,y^2,z^2)$, thus $I^{(2)}_{(Z,\o)}=(x,y)^2\cap (y,z)^2\cap (x,z)^2$.
\eex
Such pathologies do not occur when $(Z,\o)$ is a complete intersection:
\bel\label{Thm.Symbolic.Powers} If $(Z,\o)\sset(\k^N,\o)$ is a complete intersection (not necessarily reduced) then
 $I^{(m)}_{(Z,\o)}=I^m_{(Z,\o)}$ for any $m\in\Z_{>0}$.
\eel
\bpr Let $J\sset R$ be the defining ideal of $(Z,\o)\sset(\k^N,\o)$.
\bei
\item
If $J\sset R$ is prime then we can use the general proposition 3.5.12 of \cite{Vasconcelos}:
\beq
\text{if $R$ is Cohen-Macaulay and $J\sset R$ is a prime complete intersection then $J^{(m)}=J^m$ for any $m\in\Z_{>0}$.}
\eeq
\item If $J$ is not prime, but $R$ is regular and in the primary decomposition, $J=\cap \cp_i$, all the minimal primes $\cp_i$ are complete intersections,
 then one can use:
\beq
J^{(m)}=\cap \cp^{(m)}_i=\cap \cp^{m}_i=J^m.
\eeq

\item In general, the minimal primes $\cp_i$ are not complete intersections, then one argues as follows. Suppose $R$ is a regular
 local ring and $J$ is a complete intersection,
 with $\sqrt{J}\neq\cm$.
  Consider the ideal $\quotients{J^{(m)}}{J^m}\sset\quotients{R}{J^m}$. By the definition of symbolic powers, for any $\cp_i$ the localization vanishes:
 $(\quotients{J^{(m)}}{J^m})_{\cp_i}=\{0\}$. Thus $\quotients{J^{(m)}}{J^m}$ is a torsion. But, as $J\sset R$ is a complete intersection, the ring
  $\quotients{R}{J^m}$ has no torsion. Thus  $\quotients{J^{(m)}}{J^m}=0$.
\epr\eei

\subsection{The functor of associated graded modules}
Fix a (commutative, associative) ring $R$, and an ideal $I\sset R$. This ideal induces the filtration, $R=I^0\supset I\supset\cdots$.
 Take the associated graded ring,
$gr_{I}(R)=\oplusl_{j\ge0}\quotients{I^j}{I^{j+1}}$. Explicitly, fix some generators, $\{g_i\}$, of $I$, and its module of relations,
 $Syz(\ug)=\{\{a_i\}|\ \sum a_i g_i=0\in R\}$. Then
 \beq\label{Eq.graded.ring}
 gr_IR=\quotient{\quotients{R}{I}[\uy]}{\Big\{\sum a_i y_i,\ \ua\in Syz(\ug)\Big\}}.
\eeq
Thus $gr_{I}R$ is a graded algebra over $\quotients{R}{I}$,
  and  $Spec(gr_I R)$ is an affine scheme over $Spec(\quotients{R}{I})$.

\

Consider the category  $mod_{filt}(R)$, of (finitely generated) filtered $R$-modules,
\beq
\cdots\supset M_{i-1}\supset M_{i}\supset M_{i+1}\supset\cdots,\quad\quad I^jM_i\sseteq M_{i+j}.
\eeq
 The morphisms here are the filtered homomorphisms:
\beq
Hom_{filtr}(\{M_i\},\{N_i\})=\{\phi\in Hom_R(M,N)|\ \phi(M_i)\sseteq N_i\}.
\eeq
To each filtered $R$-module one associates a graded module over $gr_IR$, by $gr(M):=\oplusl_{j\ge0}\quotients{M^j}{M^{j+1}}$.
 The filtered morphisms of $mod_{filt}(R)$ are then sent to the graded morphisms of $mod_{gr}(gr_IR)$.  This defines the ``associated graded" functor
\beq
mod_{filt}(R)\stackrel{gr}{\to}mod_{gr}(gr_I(R)).
\eeq
 In our case $R$ is Noetherian and all the filtration are exhaustive ($\cup M_i=M$) and separated ($\cap M_i=\{0\}$), in particular  $\capl_j I^j=\{0\}$.
   Therefore this functor is faithful,
  i.e. $gr(M)=\{0\}$ implies $M=\{0\}$, \cite[Proposition I.4.1]{Nastasescu-Van Oystaeyen}.

   This functor is not exact, however it preserves exactness of strictly filtered sequences. In more detail, take a filtered morphism
 $M\stackrel{\phi}{\to} N$, i.e. $\phi(M_i)\sseteq N_i$. This morphism is called {\em strictly filtered} if $\phi^{-1}(N_i)=M_i$.
\bthe\label{Thm.graded.functor.exact}  \cite[Theorem I.4.4]{Nastasescu-Van Oystaeyen}
Consider a filtered sequence in $mod_{filt}-R$ and the associated sequence
  in $mod_{gr}(gr_I(R))$:
\[ \ *:\ L\to M\to N\quad\quad \rightsquigarrow\quad\quad
gr(*):\ gr(L)\to gr(M)\to gr(N).
\]
  Suppose all the filtrations are exhaustive and complete. Then $gr(*)$ is exact iff $*$ is exact and strictly filtered.
\ethe

\subsection{The normal cone}\label{Sec.Background.Normal.Cone}

Given a filtration $R=I^0\supset I^1\supset\cdots$ and an element $f\in R$, fix the order $p=ord_I(f)$, i.e. $p$ with $f\in I^p\smin I^{p+1}$.
The leading term of $f$ is defined as the residue class  $l.t.(f)\in \quotients{I^p}{I^{p+1}}\sset gr_I(R)$.
 We associate with  an ideal $J\sset R$ the ideal $gr_I(J)\sset gr_IR$, generated by the leading terms of all the elements of $J$.

In our case, for a triple of germs, $(Z,\o)\sset(X,\o)\sset Spec(R)$, we have the diagram:
\beq\bM R& \rightsquigarrow& gr_{I_{(Z,\o)}} R=&\oplusl_{j\ge0}\quotients{(I_{(Z,\o)})^j}{(I_{(Z,\o)})^{j+1}} \\
\cup&&\cup
\\I_{(X,\o)}&\rightsquigarrow&gr_{I_{(Z,\o)}} I_{(X,\o)}\eM
\eeq

One can write explicitly:
$gr_{I_{(Z,\o)}} I_{(X,\o)}=\oplusl_{j\ge0}\quotient{I_{(X,\o)}\cap(I_{(Z,\o)})^j+(I_{(Z,\o)})^{j+1}}{(I_{(Z,\o)})^{j+1}}$.

Note that the  transition $f\rightsquigarrow l.t.(f)$,  $I_{(X,\o)}\rightsquigarrow gr I_{(X,\o)}$, is not a homomorphism and is never injective/surjective.
 (Its image is the disjoint union of all the homogeneous components of $gr I_{(X,\o)}$.)

\bed
1. The ideal $gr_{I_{(Z,\o)}} I_{(X,\o)}\sset gr_{I_{(Z,\o)}} R$ is called the co-normal ideal.
\\2. The normal cone of $(X,\o)$ along $(Z,\o)$ is the scheme $\cN_{(X,\o)/(Z,\o)}=V(gr_{I_{(Z,\o)}} I_{(X,\o)})\sset Spec(gr_{I_{(Z,\o)}} R)$.
\eed

\subsubsection{Example: $(Z,\o)\sset Spec(R)$ is a complete intersection}
Let  $(Z,\o)$ be a complete intersection (not necessarily reduced) and fix some regular sequence of generators $I_{(Z,\o)}=(g_1,\dots,g_k)$.
Then the only relations among $\{g_i\}$ are the Koszul relations, therefore equation \eqref{Eq.graded.ring} gives:
\beq
gr_{I_{(Z,\o)}} R\approx \cO_{(Z,\o)}[\uy],\quad \text{ with } \uy=(y_1,\dots,y_k), \quad \text{ and }\quad Spec(gr_{I_{(Z,\o)}} R)=(Z,\o)\times \k^k.
\eeq
(Here the isomorphism is defined by the choice of the generators $\{g_i\}$.
 This ambiguity results in the action $GL_{\cO_{(Z,\o)}}(k)\circlearrowright(Z,\o)\times\k^k$, see also below.

 For any $f\in I^p_{(Z,\o)}\smin I^{p+1}_{(Z,\o)} \sset R$ we have:
\beq
f=\suml_{\sum m_j=p}g^{m_1}_1\dots g^{m_k}_k a_{m_1\dots m_k}, \ \text{ for some elements } \ \{a_{m_1\dots m_k}\in  R\}.
\eeq
Using this expansion we write down the leading term of $f$ explicitly:
\beq
\tf=\sum y^{m_I}_Ia_{m_I}|_{(Z,\o)}\in\cO_{(Z,\o)}[\uy], \text{ where } \cO_{(Z,o)}=\quotients{R}{I_{(Z,\o)}}\sset  gr_{I_{(Z,\o)}} R
\eeq
(By the construction: $\tf\neq0\in\cO_{(Z,\o)}[\uy]$, therefore $ord_{(Z,\o)}(\tf)=p$ too.)

 The coefficients $\{a_{m_1\dots m_k}\}$ are not unique, because of the Koszul relations.
But the restrictions $\{a_{m_I}|_{(Z,\o)}\in \cO_{(Z,o)}\}$ are defined uniquely.
   (Indeed, if $\sum \uy^{m_I}_Ia_{m_I}|_{(Z,\o)}=\sum \uy^{m_I}_Ib_{m_I}|_{(Z,\o)}$ then $\sum g^{m_I}_I(a_{m_I}-b_{m_I})=0\in\quotients{R}{I^{p+1}}_{(Z,\o)}$,
    which means a syzygy between $g_i$'s
    in $\quotients{R}{I^{p+1}_{(Z,\o)}}$. It lifts to a syzygy in $R$, with some contribution of term from $I^{p+1}_{(Z,\o)}$.
    But,
  $\{g_i\}$  being a regular sequence,  all the syzygies are linear combinations of the Koszul ones, and this forces $a_{m_I}-b_{m_I}$ to
  belong to $I_{(Z,\o)}$, which implies    $a_{m_I}|_{(Z,\o)}-b_{m_I}|_{(Z,\o)}=0\in\cO_{(Z,\o)}$.)

\

The set of all such leading terms, $I_{(X,\o)}\ni f\rightsquigarrow \tf:=l.t.(f)$, generates the co-normal ideal $gr_{I_{(Z,\o)}} I_{(X,\o)}\sset \cO_{(Z,\o)}[y_1,\dots,y_k]$.
 Again, the transition $I_{(X,\o)}\rightsquigarrow gr_{I_{(Z,\o)}} I_{(X,\o)}$ is not a homomorphism and is never injective/surjective.
  The ideal $gr_{I_{(Z,\o)}} I_{(X,\o)}$
  is graded (by construction).

\bex
Suppose $R=\cO_{(\k^N,o)}$ and  $(Z,\o)\sset (\k^N,o)$ is just a reduced point, then $I_{(Z,\o)}=\cm\sset \cO_{(\k^N,\o)}$ and
$gr_{(Z,\o)}\cO_{(\k^N,\o)}=\oplusl_{j\ge0} \quotients{\cm^j}{\cm^{j+1}}$.
 The transition $I_{(X,\o)}\rightsquigarrow gr{(Z,\o)}I_{(X,\o)}$ takes the leading term of Taylor expansion, $f=f_p+f_{>p}\rightsquigarrow f_p$. The normal cone is just the tangent cone,
 $T_{(X,\o)}\sset T_{(\k^N,\o)}=Spec(\k[y_1,\dots,y_N])$.
\eex

To identify $gr_{I_{(Z,\o)}} R\approx \cO_{(Z,\o)}[\uy]$  we have chosen a set of generators of $I_{(Z,\o)}$, but the
  dependence of the image of $gr_{(Z,\o)} I_{(X,\o)}$ in $\cO_{(Z,\o)}[\uy]$ on this particular choice is non-essential:
\bel\label{Thm.Normal.Cone.is.well.defined}
The ideal $gr_{(Z,\o)} I_{(X,\o)}$ in $\cO_{(Z,\o)}[\uy]$ is well defined up to the action $GL_{\cO_{(Z,o)}}(k)\circlearrowright (Z,\o)\times\k^k$, in particular
 the subscheme  $\cN_{(X,\o)/(Z,\o)}\sset (Z,\o)\times\k^k$ is defined up to an isomorphism.
\eel
\bpr
Fix some other set of generators $\{g'_i\}$ of $I_{(Z,\o)}$, then $g'_i=\suml_j \phi_{ij}g_j$, where the matrix $\{\phi_{ij}\}\in Mat_{k\times k}(R)$
 is invertible. This matrix induces an automorphism of $\cO_{(Z,\o)}[\uy]$, denote it by $\phi$.
 Then any homogeneous expansion $\tf=\sum y^{m_I}_Ia_{m_I}$ transforms to   $\tf=\sum \phi(y^{m_I}_I)a_{m_I}$ Thus
  $gr_{(Z,\o)} I_{(X,\o)}(\{g_i\})$ and    $gr_{(Z,\o)} I_{(X,\o)}(\{g'_i\})$ differ by an element of $GL_{\cO_{(Z,o)}}(k)$.
\epr

\subsection{Strictly complete intersections}\label{Sec.Background.Strictly.Complete.Intersections}
The tangent cone of a hypersurface germ is a hypersurface, but the tangent cone to a complete intersection is not necessarily a complete intersection.
\bex\label{Ex.ci.whose.PT.is.smooth.not.ci} For the
complete intersection $(X,\o)=\{x^2+zy^3=xy+z^3=0\}\sset(\k^3,\o)$ the
tangent cone is $T_{(X,\o)}=\{x^2=xy=xz^3=z^6-zy^5=0\}$. Indeed, \cite[\S15.10.3]{Eisenbud-book}, it
is enough to check
the Groebner basis of the homogenized ideal, $\{w^2x^2+zy^3,wxy+z^3\}$,
with respect to any monomial ordering. For the ordering $x>y>z>w$ the Groebner basis is:
\beq
\{w^2x^2+zy^3,wxy+z^3,wxz^3-zy^4,z^6-y^5z\}.
\eeq
By sending $w\to1$ and taking the leading terms we get the projectivized tangent cone.
Now, by direct check, this projectivization $\P T_{(X,\o)}\sset\P^2$ is
a collection of {\em smooth} (!) points, whose defining ideal is not a complete intersection.
\eex
For various other pathologies of tangent cone and conditions to prevent them see \cite{Heinzer-Kim-Ulrich}.
\bed\label{Def.SCI}
The germ $(X,\o)\sset Spec(R)$ is called a strictly complete intersection (s.c.i.) over $(Z,\o)$ if the normal cone,
 $\cN_{(X,\o)/(Z,\o)}$, is a complete intersection over $(Z,\o)$. Algebraically: the ideal
  $gr_{(I_{(Z,\o)})} I_{(X,\o)}\sset gr_{(I_{(Z,\o)})}(R)$ is generated by a regular sequence.
\eed
%Yet another name is: $I_{(X,\o)}$ is generated by a {\em super-regular} sequence, \cite{Valabrega-Valla}.

Many results around this notion are scattered in the literature. We collect here the relevant results and examples.
\bex\label{Ex.Hypersurface.not.always.SCI}
Let $(X,\o)\sset(\k^N,\o)$ be a hypersurface with $(Z,\o)\sset (X,\o)$ a  complete intersection.
\bee[i.]
\item  Suppose $(Z,\o)$ is irreducible and reduced, then $(X,\o)$ is s.c.i. over $(Z,\o)$.
 Indeed: if $I_{(X,\o)}=(f)$ then $gr_{(Z,\o)}I_{(X,\o)}=(\tf)$, and $\tf\in \cO_{(Z,\o)}[\uy]$ is regular, i.e. not a zero divisor.
\item  If $(Z,\o)$ is reduced but
 reducible then $gr_{(Z,\o)}I_{(X,\o)}$ is still generated by one element, but this might be not a regular sequence, since this element can be a zero divisor.
  More precisely, let $I_{(X,\o)}=(f)$ with $(Z,\o)=\cup (Z_i,\o)$ and $f\in \cap (I^{p_i}_{(Z_i,\o)}\smin I^{p_i+1}_{(Z_i,\o)})$. Then
 $(X,\o)$  is s.c.i. over $(Z,\o)$ iff $p_1=\dots=p_n$.
\item  Similarly for the case: $(Z,\o)$ is a multiple of an irreducible germ, e.g. $I_{(Z,\o)}=(x^2_1,x^2_2)\sset\k[[x_1,x_2,x_3]]$
 and $I_{(X,\o)}=(x^3_1+x^3_2)\sset \k[[x_1,x_2,x_3]]$.
   Here the ideal $gr_{(Z,\o)}I_{(X,\o)}=(x_1y^2_1+x_2y^2_2)\sset\cO_{(Z,\o)}[\uy]$ is principal. But its generator is not regular, being a zero divisor.
\eee\eex

\bex
Suppose $(Z,\o)$ is just a reduced point, then definition \ref{Def.SCI} reads as follows:
 the germ $(X,\o)$ is called a {\em strictly complete
intersection}, s.c.i. at $\o$, if it is a complete intersection and its
tangent cone is a complete intersection too. (The condition ``$(X,\o)$ is a complete intersection at $\o$" is redundant here, as we show below.) Thus a hypersurface germ is always a s.c.i. at the origin.
 The name ``strict complete intersection" seems to be coined
 by \cite[pg.31]{Benett-1977}. The name ``strong complete intersection" is used
 in commutative algebra to denote ``geometric" complete intersections, i.e. rings of the
 form $\quotients{S}{(f_1,\dots,f_r)}$, where $S$ is a regular local ring and $\{f_i\}$ is a regular sequence, \cite{Heitmann-Jorgensen}.
 The name ``absolute complete intersection" would suggest that both the germ and {\em all} its proper
 transforms and exceptional loci in the resolution are locally complete intersections.
\eex
\bprop\label{Thm.SCI.over.germ.implies.CI.at.a.point}
If $(X,\o)$ is a strictly complete intersection over $(Z,\o)$ then $(X,\o)$ is a complete intersection,
 i.e.  $I_{(X,\o)}$ is generated by a regular sequence. Moreover, there exists a choice of generators,  $I_{(X,\o)}=(f_1,\dots,f_r)$, such
  that the leading terms, $\{\tf_i\}$, form a regular sequence that generates  $gr_{(Z,\o)}I_{(X,\o)}$.
\eprop
(For the proof see Corollary 2.4 of \cite{Valabrega-Valla}. Following that paper the sequence $\{f_i\}$ is often called a ``super-regular" sequence.)

This proposition, together with example \ref{Ex.ci.whose.PT.is.smooth.not.ci}, show that the condition ``$(X,\o)$ is a s.c.i. over $(Z,\o)$'' is  stronger
 than the condition ``$(X,\o)$ is a complete intersection as a scheme over $(Z,\o)$''.

\bex Suppose  $(Z,\o)\sset(\k^N,\o)$ is a reduced complete intersection such that $(Z,\o)\sset(X,\o)$ is also a complete intersection,
 and $(X,\o)\sset(\k^N,\o)$ is also a complete intersection.
(Equivalently, $I_{(X,\o)}\sset \cO_{(\k^N,\o)}$ has a basis that can
be extended to a basis of $I_{(Z,\o)}$. Indeed, choose a regular sequence that generates the defining ideal of $(Z,\o)\sset(X,\o)$, take
 some representatives $g_1,\dots,g_{k-r}\in \cO_{(\k^N,\o)}$. Take some generators $f_1,\dots,f_r$ of $I_{(X,\o)}\sset \cO_{(\k^N,\o)}$, then
  $I_{(Z,\o)}\sset \cO_{(\k^N,\o)}$ is generated by $\{f_j\},\{g_i\}$. And this is a regular sequence.)
 Then $(X,\o)$ is s.c.i. over $(Z,\o)$. Indeed,
 take a basis of $I_{(Z,\o)}\sset \cO_{(\k^N,\o)}$ as above, $\{f_j\},\{g_i\}$.
  Then $\tf_1$,\dots,$\tf_r$ is a regular sequence in $\cO_{(Z,\o)}[y_1,\dots,y_k]$ that generates $gr_{(Z,\o)}I_{(X,\o)}$.
\eex

\

If  $(X,\o)$ is smooth and $(Z,\o)\sset(\k^N,\o)$ is a reduced complete intersection the the condition `` $(Z,\o)\sset(X,\o)$ is also a complete intersection"
 holds automatically. More generally, this often holds when $(Z,\o)\not\sseteq Sing(X,\o)$.
Usually we assume $(Z,\o)\sseteq Sing(X,\o)$.
\bex
Let $I_{(X,\o)}=(x^2w+y^4+z^4,y^2w+x^4)$ and $I_{(Z,\o)}=(x,y,z)$. Here $(Z,\o)=Sing(X,\o)$.
 Then $gr_{(Z,\o)}I_{(X,\o)}=(x^2w,y^2w,x^6-y^2(y^4+z^4))\sset\cO_{(\k^1,0)}[x,y,z]$ is not
a complete intersection.  Note though that $gr_{(Z,\o)}I_{(X,\o)}$ is {\em generically} complete intersection along $(Z,\o)$.
\eex
Example \ref{Ex.Hypersurface.not.always.SCI} shows that ``being s.c.i. at $\o$'' does not imply ``being s.c.i. over $(Z,\o)$''.
The converse implication does not hold either:
\bex
Fix any complete intersection $(Z,\o)\sset(\k^N,\o)$ and some basis $I_{(Z,\o)}=(g_1,\dots,g_k)$. Let $f_1(y_1,\dots,y_k)$, \dots, $f_r(y_1,\dots,y_k)$
 be some homogeneous polynomials, $r\le k$, such that $(f_1,\dots,f_r)$ is a regular sequence in $\k[y_1,\dots,y_k]$. Consider the ideal
\[
I_{(X,\o)}=(f_1(g_1,\dots,g_k),\dots,f_r(g_1,\dots,g_k)).
\]
 Then $(X,\o)$ is a complete intersection at $\o$ and s.c.i. over $(Z,\o)$. But in general
   $(X,\o)$ is not s.c.i. at $\o$. As a particular example, let $I_{(Z,\o)}=(x^2+zy^3,xy+z^3)$, see example \ref{Ex.ci.whose.PT.is.smooth.not.ci},
     and take $f_1=t_1$, $f_2=t_2$, i.e. $I_{(X,\o)}=I_{(Z,\o)}$.
    Then $gr_{(Z,\o)} I_{(X,\o)}=(y_1,y_2)\sset\cO_{(Z,\o)}[y_1,y_2]$, thus $(X,\o)$ is s.c.i. over $(Z,\o)$.
     But $(X,\o)$ is not s.c.i. at the origin.
\eex

\bex
A warning: even if $(X,\o)$ is a s.c.i. at $\o$, $(Y,\o)$ is smooth, and the intersection $(X\cap Y,\o)$ is proper, it can happen that
 $(X\cap Y,\o)$
 is in general not a s.c.i. at $\o$.
 For example, let $(X,\o)=\{x^2y+z^3+h.o.t.=0=y^3+z^3+h.o.t.\}\sset(\k^3,\o)$, where $h.o.t.$ denote some elements of $(x,y,z)^4$.
  Let $(Y,\o)=\{z=0\}$, then $(X\cap Y,\o)=\{z=x^2y+h.o.t.=0=y^3+h.o.t.\}\sset(\k^3,\o)$ is not a s.c.i.
\eex

\

\subsection{Good bases and multiplicity sequences}\label{Sec.Background.Good.Bases}

Let $(Z,\o)$ be a reduced complete intersection and
$(X,\o)$ be s.c.i. over $(Z,\o)$. By  Proposition \ref{Thm.SCI.over.germ.implies.CI.at.a.point} we can choose some regular sequence of generators,
 $f_1,\dots,f_r\in I_{(X,\o)}$, whose leading terms form a (graded) basis $\{\tf_i\}$ of $gr_{(Z,\o)}I_{(X,\o)}$.
   We can assume in addition:
 \beq\label{Eq.Good.Basis.orders}
 ord_{(Z,\o)}(f_1)\le ord_{(Z,\o)}(f_2)\le\cdots\le ord_{(Z,\o)}(f_r).
 \eeq
By the construction: $ord_{(Z,\o)}(f_i)=ord_{(Z,\o)}(\tf_i)$.
Recall that for complete intersections the ordinary and symbolic powers coincide, $I^m_{(Z,\o)}=I^{(m)}_{(Z,\o)}$.
Then $ord_{(Z,\o)}(f_i)$ is the generic order of vanishing of $f_i$ along $(Z,\o)$, in particular it is the same on all the components of $(Z,\o)$.
 Indeed, (cf. example \ref{Ex.Hypersurface.not.always.SCI}) let $(Z,\o)=\cup(Z_j,\o)$ be the irreducible decomposition,
  so that $f_i\in \capl_j I^{(ord(f_i))}_{(Z_j,\o)}$.
  If $f_i\in I^{(ord(f_i)+1)}_{(Z_j,\o)}$   holds for some $j$, then $\tf_i\in \quotients{I^{(ord(f_i))}_{(Z,\o)}}{I^{(ord(f_i)+1)}_{(Z,\o)}}=
   \quotients{I^{ord(f_i)}_{(Z,\o)}}{I^{ord(f_i)+1}_{(Z,\o)}}\sset gr_{(Z,\o)}(R)$ is a zero divisor,
   contradicting regularity of the sequence $\tf_1$,\dots, $\tf_r$.
\bed
The sequence $f_1,\dots,f_r$, as in equation \eqref{Eq.Good.Basis.orders}, is called {\em a good basis} of $I_{(X,\o)}$.
 The sequence of integers
  $(ord_{(Z,\o)}(f_1)$, $ord_{(Z,\o)}(f_2)$,\dots,$ord_{(Z,\o)}(f_r))$ is called {\em the multiplicity sequence}  of $(X,\o)$ along $(Z,\o)$,
   associated with $\{f_i\}$.
\eed
\hspace{-0.4cm}\parbox{12.7cm}{\bex
(The case: $(Z,\o)$ is a point.) Let $(X,\o)\sset(\k^{n+r},\o)$ be a s.c.i. with multiplicity sequence $(p_1,\dots,p_r)$. Blowup at $\o$, see the diagram, then
 $\tX\cap E\sset E\approx\P^{n+r-1}$ is a globally complete intersection of multi-degree  $(p_1,\dots,p_r)$.
\eex}\quad\quad$\bM \tX\cap E\sset\tX\sset Bl_\o(\k^{n+r},\o)\\ \downarrow\quad\quad\quad\downarrow\quad\quad\quad \ \
\downarrow\quad\quad\\\o \ \ \in(X,\o)\sset(\k^{n+r},\o)\eM$.

Eventhough the good basis is never unique, the multiplicity sequence is well defined.
\bprop
1. The multiplicity sequence of $(X,\o)$ along $(Z,\o)$ does not depend on the choice of bases of $I_{(X,\o)}$, $\tI_{(X,\o)}$.
\\2.  The product $\prod ord(f_i)$ equals the generic multiplicity of $(X,\o)$ along $(Z,\o)$. In particular, it is the same on all the components of $(Z,\o)$.
\eprop
(For the proof of part, for the case of point, $(Z,\o)=\o$, see example 12.4.9 in \cite{Fulton-book}.)
\bpr
We decompose $I_{(Z,\o)}=\cap\cp_i$ and localize at $\cp_i$. We get the local ring $(\cO_{(\k^{n+r},\o)})_{\cp_i}$
with the maximal ideal $(I_{(Z_i,\o)})_{\cp_i}$ and elements $\big\{f_j\in (I_{(Z_i,\o)})^{ord(f_j)}_{\cp_i}\big\}$.
 So the situation is reduced to the case when   $(Z,\o)$ is a point.

Suppose $(Z,\o)$ is a point, then $\{\tf_i\}$ is the graded basis  of the tangent cone,  $I_{T_{(X,\o)}}$.
 Thus the uniqueness of the multiplicity sequence (up to permutation) follows from the fact that any two
 graded bases of the graded ideal $I_{T_{(X,\o)}}$ are related by a graded(!) invertible linear map.

To compute the multiplicity of $(X,\o)$, at $\o$, let $l_1,\dots,l_n\in\cm\sset \cO_{(\k^{n+r},\o)}$ be some generic elements.
 Geometrically they define a smooth subspace transversal and complementary to $(X,\o)$. The quotient ring
  $\quotients{\cO_{(\k^{n+r},\o)}}{(l_1,\dots,l_n,f_1,\dots,f_r)}$ is still a strictly complete intersection. Therefore we can assume from
   the beginning $n=0$, i.e. $(X,\o)$ is a one-point scheme and $mult(X,\o)=dim\quotients{\cO_{(\k^{r},\o)}}{(f_1,\dots,f_r)}$. But then
    the multiplicity can be computed by blowing up. And the total transform of $\o$ is a complete intersection of multidegree $m_1,\dots,m_r$, thus of
     degree $\prodl^r_{i=1}m_i$.
\epr

\bex (Behavior of multiplicity sequence in a family.)
Let $(X,\o)$ be a s.c.i. over $(Z,\o)$ with a good basis $\{f_1,\dots,f_r\}$. Consider a
deformation $\{f_1+\ep g_1,\dots,f_r+\ep g_r\}_{\ep\in(\k^1,0)}$ that preserves the multiplicity sequence  (generic multiplicity over $(Z,\o)$),
 $mult(g_i)\ge mult(f_i)$. Then, by the openness of regularity in deformations, the generic member of this family is a s.c.i. over $(Z,\o)$.
  (Indeed, any relation $\sum r_i(f_i+\ep g_i)=0$ leads  to the relation $\sum r_i|_{\ep=0}f_i=0$, which is necessarily Koszul.
   Subtract this relation from the initial one, to get $r_i|_{\ep=0}=0$, i.e. $r_i$ is divisible by $\ep$. Divide all $r_i$ by (a power of $\ep$) and
    proceed in the same way. The statement then follows by Nakayama-type argument.)

\

If the multiplicities are not preserved then a flat deformation of s.c.i. is not s.c.i. For example,
 the family of ideals $I=(x^3,tx^2+y^4)\sset \k[[x,y,t]]$ defines
 s.c.i. at the origin for $t=0$ but not for $t\neq0$.
\eex

\subsection{Singularities generically ordinary along $(Z,\o)$}\label{Sec.Background.Ordinary.Singularities} Recall that an isolated hypersurface
singularity, $\{f_p+f_{>p}=0\}\sset(\k^{n+1},\o)$, is called {\em an \omp} if its projectivized
tangent cone, $\{f_p=0\}\sset\P^n$, is smooth.
In the case $\k=\C$ this can be stated as follows: the hypersurface germ is topologically equisingular to
 $\{\suml_i x^p_i=0\}\sset(\C^{n+1},\o)$.
 This is ``the lowest/simplest'' hypersurface singularity of a given multiplicity. Similarly, among the s.c.i. germs of a given multiplicity at $\o$,
  the ``lowest'' is the one whose projectivized tangent cone is a smooth complete intersection.

\

Let $(Z,\o)$ be reduced, $dim(Z,\o)>0$, though $(Z,\o)$ is not necessarily a complete intersection or pure dimensional.
 Then the primary decomposition contains only prime ideals, $I_{(Z,\o)}=\cap \cp_i$.
\bed
$(X,\o)$ is called generically ordinary along $(Z,\o)$ if for any $i$ the localization $(gr_{(Z,\o)}I_{(X,\o)})_{(\cp_i)}$ is a complete intersection over
 $Spec((\cO_{(Z,\o)})_{(\cp_i)})$,  whose projectivization is smooth.
\eed
\bex
If $(Z,\o)\not\sseteq Sing(X,\o)$ and is irreducible then $(X,\o)$ is generically smooth along $(Z,\o)$, in particular $(X,\o)$
 is generically ordinary along $(Z,\o)$. However, in this case $\{ord_{(Z,\o)}(f_i)=1\}$, thus, as will be
  shown in example \ref{Ex.Discriminant.Classical.is.empty.for.deg=1}, the discriminant is empty.
   Therefore we assume $(Z,\o)\sseteq Sing(X,\o)$.
\eex
\bex
$\bullet$ Let $I_{(X,\o)}=(x^pz+y^{p+1})$, $p\ge2$, and $I_{(Z,\o)}=(x,y)$. Then $gr_{(Z,\o)}I_{(X,\o)}=(x^pz)\sset\cO_{(\k^1,\o)}[x,y]$ and
  $(gr_{(Z,\o)}I_{(X,\o)})_{(x,y)}=(x^p)\sset(\cO_{(\k^1,\o)}[x,y])_{(x,y)}$. Thus $(X,\o)$ is not generically ordinary along $(Z,\o)$.
\bei
\item   The hypersurface $\{x^pz^q+y^p+x^{p+1}=0\}$ is generically ordinary along $(Z,\o)=\{x=0=y\}$.
\item
Let $I_{(X,\o)}=(x^2z+y^2w+x^4+y^4)$ and $I_{(Z,\o)}=(x,y)$. Then $gr_{(Z,\o)}I_{(X,\o)}=(x^2z+y^2w)$ and
\[
(gr_{(Z,\o)}I_{(X,\o)})_{(x,y)}=(x^2z+y^2w)\sset(\cO_{(\k^2,\o)}[x,y]))_{(x,y)}.
\] This is a hypersurface and its projectivization is smooth.
   Thus $(X,\o)$ is generically ordinary along $(Z,\o)$.
\eei
\eex
\bex
$\bullet$ Let $I_{(X,\o)}=(x^2w+y^4+z^4,y^2w+x^4)$ thus $I_{(Z,\o)}=(x,y,z)$. Then $gr_{(Z,\o)}I_{(X,\o)}=(x^2w,y^2w,x^6-y^2(y^4+z^4))$   and
$(gr_{(Z,\o)}I_{(X,\o)})_{(x,y,z)}=(x^2,y^2)\sset(\cO_{(\k^1,\o)}[x,y,z]))_{(x,y,z)}$. This is a complete intersection, but its projectivization is not a
 smooth subscheme of $\P^2_{\K}$, here $\K$ is the fraction field of $\cO_{(\k^1,\o)}$.

$\bullet$
Let $I_{(X,\o)}=((x^2+y^2+z^2)w+h_1(x,y,z),xyw+h_2(x,y,z))$, here the orders of $h_1,h_2$ are$\ge4$, and $I_{(Z,\o)}=(x,y,z)$.
 Then
\begin{multline} gr_{(Z,\o)}I_{(X,\o)}=((x^2+y^2+z^2)w,xyw,xyh_1-(x^2+y^2+z^2)h_2),\\ \text{  and } \quad
 (gr_{(Z,\o)}I_{(X,\o)})_{(x,y,z)}=(xy,x^2+y^2+z^2)\sset(\cO_{(\k^1,\o)}[x,y,z]))_{(x,y,z)}.
 \end{multline}
  This is a complete intersection and its projectivization is a
 smooth subscheme of $\P^2_{\K}$. Thus $(X,\o)$ is generically ordinary over $(Z,\o)$. Note that $(X,\o)$ is not s.c.i. at the origin.
\eex

\section{The classical discriminant of projective complete intersections}\label{Sec.Discriminant.Classical}
While there are many extensive treatments of the discriminant of projective hypersurfaces in $\P^N$,
 see e.g.  \cite{Gelfand-Kapranov-Zelevinsky},
we do not know any textbook or lecture notes on the discriminant of projective complete
intersections.

However, for some recent particular results on the classical discriminant of projective complete intersections see  \cite{Esterov},  \cite{Benoist}, \cite{C.C.D.R.S.2011}.
In particular, in many cases the multi-degrees were computed.

In \cite{Teissier1976}, \cite{Looijenga-book} one treats mostly the local case. See also \cite[\S I.2.2]{AGLV2}
 for a collection of known local facts.

 In this section we (re)prove some of the standard needed results.

\subsection{The critical locus and the discriminant of a map}\label{Sec.Background.Critical.Locus.Discrim}

 Let $X\proj S$ be a flat map  of (algebraic/analytic/formal) spaces, with fibres of pure dimension $d$.
  Then $Crit(\pi)$ is defined, see e.g. \cite[pg.587]{Teissier1976}, by the (coherent)
sheaf of ideals  $I_{Crit(\pi)}:=Fitt_d(\Om^1_{X/S})$. Here $\Om^1_{X/S}=\quotient{\Om^1_X}{\pi^*\Om^1_S}$ is the sheaf of relative differentials,
 while $Fitt_d(...)$ is the $d$'th Fitting ideal of an $\cO_X$-module, \cite[\S20]{Eisenbud-book}.

Suppose the restriction $\pi|_{Crit(\pi)}$ is finite. The discriminant of $\pi$ is defined as the image, $\De_\pi:=\pi(Crit(\pi))\sset S$,
 with the Fitting scheme structure, $I_{\De_\pi}:=Fitt_0(\pi_*\cO_{Crit(\pi)})$,
 \cite[pg. 588]{Teissier1976}.
Here $\pi_*\cO_{Crit(\pi)}$ is the pushforward of the $\cO_X$-module $\cO_{Crit(\pi)}$, while $Fitt_0(...)$ is the
minimal Fitting ideal of a module, as an $\cO_S$ module,  i.e. the ideal of maximal minors of a presentation matrix of the module.

\subsection{Assumptions on the base of the family}\label{Sec.Discriminant.Assumptions.on.the.base.family}
Consider the family of complete intersections,
\[
\cX=\{F_1(x,s)=F_2(x,s)=\dots=F_r(x,s)=0\}\sset \P^{n+r}\times S.
\]
 We denote the fiber in $\cX$ over the point $s\in S$ by $X_s\sset\P^{n+r}$.
 We have  the natural projection $\cX\stackrel{\pi}{\to}S$.
%Depending on the context, one uses different spaces to parameterize complete intersections.
%When possible we will  not specify the parameter space and denote it just by $S$.

 We assume:
 \bei
 \item
 $S$ is quasi-projective, smooth, connected.
\item The generic fibre over $S$ is a smooth complete intersection in $\P^{n+r}$ and  the  family $\cX\sset S\times\P^{n+r}$ is smooth.
\item  Denote by $\De\sset S$ the subset of points whose fibres are singular or not of expected dimension. Denote by  $\De_{A_1}\sseteq \De\sset S$ the
 subset of points corresponding to fibres with just one node. Then we assume that $\De_{A_1}$ is dense in $\De$ and is  connected in Zariski topology.
\eei
For hypersurfaces of degree $p$ in $\P^{n+1}$ the standard parameter space is $|\cO_{\P^{n+1}}(p)|$. For complete intersections
of  multi-degree $(p_1,\dots,p_r)$ in $\P^{n+r}$ one can consider the multi-projective space $\prodl_i|
\cO_{\P^{n+r}}(p_i)|$. To a point of this space, $(f_1,\dots,f_r)\in \prodl_i |\cO_{\P^{n+r}}(p_i)|$,
 corresponds the subscheme $X_{\uf}=\{f_1=\dots=f_r=0\}\sset\P^{n+r}$. These subschemes are projective
complete intersections when the polynomials $\{f_i\}$ form a regular sequence. Thus
there exists a Zariski open subset $\cU\sset\prodl_i|\cO_{\P^{n+r}}(p_i)|$,
whose points correspond to projective complete intersections.
(Note that the complement, $\prodl_i|\cO_{\P^{n+r}}(p_i)|\smin \cU$,  is of high codimension.)
 This is the reason to consider $S$ as a
{\em parameter space} for globally complete intersections, eventhough for $r>1$ the correspondence $S\ni\uf
\rightsquigarrow X_{\uf}\sset\P^{n+r}$ is far from being injective.

\subsection{An example: the critical locus in the case $S=\prodl_i| \cO_{\P^{n+r}}(p_i)|$}\label{Sec.Discriminant.Classical.Critical.Locus.Explicit}

 The critical locus of $\pi$ is defined
 (as in \S\ref{Sec.Background.Critical.Locus.Discrim})
 by the sheaf of ideals $Fitt_n(\Om^1_{\cX/S})\sset\cO_\cX$.
We write down the generators $Fitt_n(\Om^1_{\cX/S})$ explicitly.

 Fix some points $x\in \P^{n+r}$, $s\in S$
 and work locally near these points, with the local coordinates $x=(x_1,.\dots,x_{n+r})$, $s=(s_1,\dots,s_r)$. We work with modules and then glue them to sheaves.
 One has
 \[\Om^1_{(\cX,(x,s))}=\quotient{\cO_{(\cX,(x,s))}(dx_1,\dots,dx_{n+r},d s_1,\dots,d s_r)}{dF},\quad
   \pi^*\Om^1_{(S,s)}=\cO_{(\cX,(x,s))}(d\us),
 \]
 where the differentials in $dF$ are taken with respect to both variables, $(x,s)$. Therefore
 \beq
 \Om^1_{(\cX/S,(x,s))}=\quotient{\cO_{(\cX,(x,s))}(dx_1,\dots,dx_{n+r})}{\{\suml_i \frac{\di F_j}{\di x_i} dx_i\}_{j=1,\dots,r}}
 \eeq
The $\cO_{(\cX,(x,s))}$-resolution of this module begins as
\beq\label{Eq.Relative.Differentials}
\to \cO^{\oplus r}_{(\cX,(x,s))}\stackrel{\{\di_i F_j\}_{ij}}{\to}\cO^{\oplus(n+r)}_{(\cX,(x,s))}\to \Om^1_{\cX/S}\to0.
\eeq
Therefore the ideal $Fitt_n(\Om^1_{\cX/S})\sset \cO_{(\cX,(x,s))}$ is defined by all the $r\times r$-minors of the matrix
\beq\label{Eq.Critical.Locus.Local}
\bpm \di_1 F_1&\dots& \di_{n+r}F_1\\\dots&\dots&\dots\\\di_1 F_r&\dots&\di_{n+r}F_r\epm.
\eeq
To get the sheaf of ideals $I_{Crit(\pi)}\sset\cO_{\cX}$ we pass from the local coordinates of $\P^{n+r}$ to the homogeneous coordinates
 $[x_0:\cdots:x_{n+r}]\in\P^{n+r}$. Using Euler's formula
$\suml^{n+r}_{i=0}x_i\di_i F_j=p_j F_j$ we get on $\cX$: the rows of the matrix in equation \eqref{Eq.Critical.Locus.Local} are linearly dependent iff the  extended rows of derivatives
in homogeneous coordinates, $\di_0 F_j\dots \di_{n+r}F_j$,  are linearly dependent. Therefore the explicit equations of the critical locus are:
\beq\label{Eq.Critical.Locus.Global}
Crit(\pi)=\Big\{F_1(x)=\cdots=F_r(x)=0,\ \ rank\bpm \di_0 F_1&\dots &\di_{n+r}F_1\\\dots&\dots&\dots\\\di_0 F_r&\dots&\di_{n+r}F_r\epm<r\Big\}\sset \P^{n+r}\times S.
\eeq

\subsection{The discriminant as the pushforward of the critical locus}\label{Sec.Discriminant.Pusforward}
Usually the projection $Crit(\pi)\to S$ is not flat over its image. More precisely, it is generically finite over its image (with varying degrees of fibres)
but the fibres over some points can be of positive dimension (when $S$ contains points whose fibres $X_s$ have
  non-isolated singularities or are not of expected dimension).
Yet this projection is proper everywhere, thus the pushforward $\pi_*(\cO_{Crit(\pi)})$ is a coherent sheaf of $\cO_S$ modules.

We work in the assumptions of \S\ref{Sec.Discriminant.Assumptions.on.the.base.family}.
\bed\label{Def.Classical.Discriminant}
The (classical) discriminant  $\De\sset S$ of complete intersections is the closure of a (algebraic) subscheme, which is defined by  the zero Fitting ideal,
 $Fitt_0(\pi_*(\cO_{Crit(\pi)}))$ at points where $\pi$ is finite.
\eed
\bex\label{Ex.Discriminant.Classical.is.empty.for.deg=1}
Suppose the multi-degree is $p_1=\cdots=p_r=1$. Then every fiber of $\cX\to S$ is smooth, in particular $\De_{A_1}=\empty$.
 Therefore, according to  our definition, $\De=\empty$.
\eex
\bprop\label{Thm.Classical.Discriminant.Properties}
1. (Set theoretically) A point $s\in S$ belongs to $\De$ iff the subscheme $X_s\sset\P^{n+r}$ is singular  or not of expected dimension.
\\2. $\De\sset S$ is a reduced irreducible Cartier divisor. The germ  $(\De,s)$ is smooth   iff the fibre $X_s$ has just one singularity of type $A_1$.
\eprop
\bpr
{\bf 1.} It is enough to check only the points of $\De$ over which $\pi$ is finite. (Indeed, the fibre over a point of $\De$ added by the closure procedure
 is the limit of  singular fibres, hence cannot be a smooth variety of expected dimension.)

For the points where $\pi$ is finite, it is enough to check the support of the module $\Om^1_{\cX/S}$. Note that the presentation of equation \eqref{Eq.Relative.Differentials}
 holds locally for any $S$, and for a fixed $s\in S$ the minors of  \eqref{Eq.Critical.Locus.Local} define  the singular set of $X_s$.
 This proves the statement.

{\bf 2.} Recall that $S$ is smooth and $\De$ is defined as the closure of a scheme well-defined on an open set above which $\pi$ is finite. Thus, to establish that $\De$ is
 a reduced Cartier  divisor, it is enough to check only those points, where $\pi$ is finite.

We should prove that the defining ideal $I_{\De}=Fitt_0(\pi_*(\cO_{Crit(\pi)}))\sset\cO_S$ is locally principal at each point.

 Note that $\cX$ is smooth and $Crit(\pi)\sset\cX$ is a determinantal subscheme of expected dimension. Therefore  $Crit(\pi)$
  is a Cohen Macaulay subscheme, \cite[Theorem 18.18]{Eisenbud-book}. As $(S,s)$
  is smooth, the module $\pi_*(\cO_{(Crit(\pi))})$ has a finite projective dimension and we use the Auslander-Buchsbaum formula,
  \cite[theorem 19.9]{Eisenbud-book}:
\beq
proj.dim(M)+depth(M)=dim(R).
\eeq
Since $\pi$ is finite, $M=\pi_*(\cO_{(Crit(\pi))})$ is a Cohen-Macaulay module over $\cO_{(S,s)}$, i.e. $depth(M)=dim(R)-1$.
 Therefore the minimal resolution of $\pi_*(\cO_{(Crit(\pi))})$ is of length one,
 \beq
0\to \cO_{(S,s)}^{\oplus N}\to \cO_{(S,s)}^{\oplus N}\to \pi_*(\cO_{(Crit(\pi))})\to0.
\eeq
Thus the presentation matrix is square and its Fitting ideal is principal.

As $\De\sset S$      is Cartier, to prove reducedness it is enough to find just one reduced point.
Suppose $X_{s_\o}$ has $A_1$-singularity at a point $x_\o\in \P^{n+r}$, here $s_\o\in S$. Then in
 some local coordinates the defining equations of $(\cX,(x_\o,s_\o))$ are: $F_1=x_1$,\dots,$F_{r-1}=x_{r-1}$, $F_r=\suml^{n+r}_{i=r}x^2_i+g(s)$, where $g(s_\o)=0$.
 As $(\cX,(x_\o,s_\o))$ is smooth, $s_\o$ is not a critical point of $g(s)$. Thus, we can choose $g(s)$ as one of the local coordinates, denote it $s_1$.

  Then using \eqref{Eq.Relative.Differentials}
  we get: $Fitt_n(\Om^1_{\cX/S})=(x_r,\dots,x_n)\sset\cO_{(\cX,(x_\o,s_\o))}$, i.e.  in the chosen coordinates
    $Crit(\pi)=V(x_1,\dots,x_n,s_1)\sset(\k^{n+r},x_\o)\times (S,s_\o)$. Therefore $\cO_{Crit(\pi)}\approx\quotients{\k[[s_1,\dots,s_r]]}{(s_1)}$ and
     $Fitt_0(\pi_*\cO_{Crit(\pi)})=(s_1)$. This ideal defines a reduced, smooth germ $(\De,s)\sset(S,s)$.

Suppose $X_s$ has a singularity other than $A_1$ then the local length of $\cO_{Crit(\pi)}$ at this point is at least two. Thus the germ $(\De,s)$ is singular.

Suppose the singular points of $X_s$ are $\{p_i\}$, then  $Fitt_0(\pi_*\cO_{Crit(\pi)})=\prod Fitt_0(\pi_*\cO_{Crit(\pi),p_i})$,
 hence the local multiplicity of $(\De,s)$ is at least the number of these points.
 Thus, if $X_s$ has more than one singular point then $(\De,s)$ is singular.

\

Finally, we know that $\De$ is a reduced divisor, thus to prove the irreducibility it is enough to check that the space $\De\smin Sing(\De)$ is connected.
 But   $\De\smin Sing(\De)=\De_{A_1}$, the open set of points corresponding to complete intersections with just one node
 and $\De_{A_1}$ is connected.
\epr

\subsection{Discriminant as a dual variety}
For a fixed tuple $p_1,\dots,p_r$ consider the multi-Veronese embedding,
\beq
\P^{n+r}\stackrel{(\nu_1,\dots,\nu_r)}{\hookrightarrow}\prodl_{j=1}^r|\cO_{\P^{n+r}}(p_j)|^\vee,
\quad
\nu_j(x_0,\dots,x_{n+r})=\{x^{m^{(j)}_0}_0\cdots x^{m^{(j)}_{n+r}}_{n+r}\}_{\substack{\sum m^{(j)}_i=p_j\\m^{(j)}_i\ge0}}\in |\cO_{\P^{n+r}}(p_j)|^\vee.
\eeq
A hyperplane in $|\cO_{\P^{n+r}}(p_j)|^\vee$ corresponds to a choice of coefficients $\{a^{(j)}_{m^{(j)}_I}\}_{\substack{\sum m^{(j)}_i=p_j\\m_i\ge0}}$,
 (up to $\k^*$-action),
 i.e. to a choice of the hypersurface $\{f_j(\ux)=0\}\sset\P^{n+r}$.  Pullback this hyperplane under the projection
 $\prod|\cO_{\P^{n+r}}(p_j)|^\vee\to |\cO_{\P^{n+r}}(p_j)|^\vee$ and denote the resulting hyperplane by $L_j$. Thus we have $r$ hyperplanes
  and the intersection $(L_1\cap \dots\cap L_r)\cap(\nu_1,\dots,\nu_r)(\P^{n+r})$ defines the subscheme $\cap\{f_j(\ux)=0\}\sset\P^{n+r}$.

  This subscheme is a smooth complete intersection (of codimension $r$) iff the intersection is transversal.
   Thus $(f_1,\dots,f_r)\in \prod|\cO_{\P^{n+r}}(p_j)|^\vee$ belongs to the discriminant iff $\capl_j L_j$ is either tangent to $(\nu_1,\dots,\nu_r)(\P^{n+r})$
    or intersects it non-properly, i.e. the resulting codimension is smaller than $r$. Thus $\De$ is the dual variety of the embedding $(\nu_1,\dots,\nu_r)(\P^{n+r})$. In particular it is a hypersurface,
     i.e. a Cartier divisor.
\

To relate this definition to the definition in \S\ref{Sec.Discriminant.Pusforward} we note that for a regular sequence  $(f_1,\dots,f_r)$ the deformation $(f_1(\ep),\dots,f_r(\ep))$
 is flat iff each $f_j(\ep)$ is flat.
And any tuple can be deformed to a tuple $(f_1(\ep),\dots,f_r(\ep))$
defining a complete intersection with isolated singularities. Therefore  the full projective discriminant
 can be obtained as the (Zariski) closure:
\beq
\De=\overline{\Big\{\ber\text{tuples $(f_1,\dots,f_r)$, defining complete intersections}\\\text{of expected codimension,
with isolated singularities}\eer\Big\}}\sset\prodl_i|\cO_{\P^N}(p_i)|.
\eeq
In particular it is irreducible and reduced and coincides with the discriminant of \S\ref{Sec.Discriminant.Pusforward}.

\subsection{The transversal multiplicity of the discriminant}\label{Sec.Discriminant.Classical.Tranversal.Multiplicity}
 Given a complete intersection germ $(X,\o)\sset(\k^{n+r},\o)$ (with isolated singularity), choose
 some {\em generic} basis $f_1,\dots,f_r\in I_{(X,\o)}$. Then, besides the ordinary Milnor number, $\mu(X,\o)$, we define the auxiliary number
  $\mu'(X,\o)=\mu(f_2=0=\dots=f_r)$. (For $r=1$ we put $\mu'=0$.)
   This $\mu'$ is well defined and depends on $(X,\o)$ only. (By the genericity of the basis one could omit any $f_i$ instead of $f_1$. )
Accordingly, for any (global) variety $X$ with ICIS, we define the total Milnor/auxiliary numbers,
\beq
\mu_{total}=\suml_{\o\in Sing(X)}\mu(X,\o),\quad\quad \mu'_{total}=\suml_{\o\in Sing(X)}\mu'(X,\o).
\eeq
\bprop\label{Thm.Classical.Discriminant.Multiplicity} Suppose
 the projection $(\cX,(x,s))\supset Crit(\pi)\to\De\sset (S,s)$ is finite at $s\in \De$ and
 $(S,s)$ contains the (germ of) miniversal deformation of $X_s$.
 Then $mult(\De,s)=\mu_{total}(X_s)+\mu'_{total}(X_s)$.
\eprop
\bpr Fix a smooth curve germ, $(C,s)\sset(S,s)$, whose tangent line does not belong to the tangent cone $T_{(\De,s)}$, then
  $mult(\De,s)=deg\big((C,s)\cap (\De,s)\big)$.
 The later degree is computed by restriction of $(\De,s)$ onto $(C,s)$. By the base-change  properties
 of Fitting ideals (i.e. the right exactness of tensor product)
  we have
\beq
\Big(Fitt_0\big(\pi_* \cO_{Crit(\pi)}\big)\Big)|_{(C,0)}=Fitt_0\big(\pi_*( \cO_{Crit(\pi)})|_{(C,0)}\big).
\eeq
 Furthermore,
 $\pi_*(\cO_{Crit(\pi)})|_{(C,0)}=\pi_*( \cO_{Crit(\pi)}|_{\pi^{-1}(C,0)})$. Therefore we can assume $(S,s)$ a smooth curve-germ.
 The module $\pi_*(\cO_{Crit(\pi)})$ is then a skyscraper at $s$ and the degree of the ideal $Fitt_0\big(\pi_* \cO_{Crit(\pi)}\big)$
   equals the length of the module $\pi_*(\cO_{Crit(\pi)})$. Finally,
\beq
length(\pi_*(\cO_{Crit(\pi)}))=h^0(\cO_{Crit(\pi)},Crit(\pi))=\suml_{pt\in  Sing(X_s)} h^0(\cO_{Crit(\pi)},pt).
\eeq
To compute $h^0(\cO_{Crit(\pi)},pt)$ we write the local presentation:
\beq
\cO_{(Crit(\pi),pt)}=\quotient{\cO_{(S,s)}[[\ux]]}{\Big(\{F_i(\ux,s)+a_i s\}_{i=1,\dots,r},\ Fitt_0\bpm dF_1\\\dots\\dF_r\epm\Big)}.
\eeq
 Here $s$ is the local parameter
  of the curve $(S,s)$, the constants $a_i\in\k$ are generic because the curve $(C,s)$ is not tangent to the discriminant.
 The differentials $dF_i$ are taken with respect to $\ux$-variables only.

  Write $F_1(\ux,s)=h_1(\ux)+s(a_1+G(\ux,s))$, with $G(\ux,s)\in(\ux)$. Then, by Nakayama lemma, we can eliminate $s$ and
   write:
\beq
\cO_{(Crit(\pi),pt)}=\quotient{\k[[x]]}{\Big( h_2,\dots, h_r,\ Fitt_0\bpm dh_1\\\dots\\dh_r\epm\Big)}.
\eeq
Here $h_1,\dots,h_r$ are some generators of $I_{(X,\o)}\sset \k[[x]$. They are generic, as the constants $\{a_i\}$ are generic.
Finally we use the L\^{e}-Greuel formula, \cite{Le},\cite{Greuel}:
\beq
\mu(h_1,\dots,h_r)+\mu(h_2,\dots,h_r)=dim\quotient{\k[[x]]}{\Big(h_2,\dots, h_r,\ Fitt_0\bpm dh_1\\\dots\\dh_r\epm\Big)}.
\eeq
Therefore
\beq
\suml_{pt\in  Sing(X_s)} h^0(\cO_{Crit(\pi)},pt)
=\suml_{pt\in Sing(X_s)} (\mu(X_s,pt)+\mu'(X_s,pt))=\mu_{total}(X_s)+\mu'_{total}(X_s).\quad
\text{\epr}\eeq

\beR\label{Thm.Classical.Discriminant.Multiplicity.one.dimensional.base}
In the last proposition $(S,s)$ was assumed ``large enough". Often the deformation space is rather small then
 the statement should be corrected. Consider a particular case, $(S,s)$ being just one dimensional, with
 $\cX\sset(\k^{n+r},\o)\times(S,s)$ defined by $I_{\cX}=\{f_1(\ux),\dots,f_{j-1}(\ux),f_j(\ux)+s,f_{j+1}(\ux),\dots,f_r(\ux)\}$.
  Here we do not assume $\{f_i\}$ to be generic, but we assume that both $\{f_i\}_{1\dots r}$ and  $\{f_i\}_{j\neq i}$
   define isolated (complete intersection) singularities.
In this case, instead of the invariant $\mu'(X,\o)$, we define the invariant
\beq
\mu_\hj:=\mu\Big(f_1(\ux),\dots,f_{j-1}(\ux),f_{j+1}(\ux),\dots,f_r(\ux)\Big).
\eeq
Then the same proof of proposition \ref{Thm.Classical.Discriminant.Multiplicity} gives:  $mult(\De,s)=\mu(X_s,\o)+\mu_\hj(X_s,\o)$.
This formula is well known, see e.g. \cite[page 589]{Teissier1976}.
\eeR

\section{The discriminant of transversal singularity type, $\Db$}\label{Sec.Db.Definition}

\subsection{The definition of  $\Db\sset Z$ as a pullback of the classical discriminant}\label{Sec.Db.Definition.via.Classical}

\subsubsection{The local case}
Let $(X,\o)\sset(\k^{n+r},\o)$ with $Sing(X,\o)=(Z,\o)$  a reduced complete intersection. Suppose $(X,\o)$ is s.c.i. over $(Z,\o)$.
 Fix a good basis, $I_{(X,\o)}=\bl f_1,\dots,f_r\br$, such that the leading terms of $(f_i)$ form a basis of $gr_{(Z,\o)} I_{(X,\o)}$,
 as in \S\ref{Sec.Background.Good.Bases}.
Let the generic multiplicity of $f_i$ along $(Z,\o)$ be $p_i$.
Fix some basis $I_{(Z,\o)}=\bl g_1,\dots,g_k\br$.
Projectivize the normal cone  to get the family:
\beq\label{Eq.Fibration.E*X->Z}
\P(\cN_{(X,\o)/(Z,\o)})=\capl_{j=1}^r\ \{\suml_{m_1+\cdots+m_k=p_j}
y^{m_1}_1\cdots y^{m_k}_ka^{(j)}_{m_1,\dots, m_k}=0\}\sset
(Z,\o)\times Proj(\k[y_1,\dots,y_k]).
\eeq

Let  $\De$ be the  {\em classical
discriminant} in the parameter space of projective complete intersections in $\P^{k-1}$ of codimension
$r$ and multidegree   $(p_1,\ldots, p_r)$, see \S\ref{Sec.Discriminant.Classical}.
It is a hypersurface, defined by one equation, $\De=\{\frD=0\}$.
 We assume that $(X,\o)$ is generically ordinary along $(Z,\o)$, see \S\ref{Sec.Background.Ordinary.Singularities}, thus
$\frD(\{a^{(j)}_{m_1,\dots,m_k}\}_{\substack{\suml_i m_i=p_j\\1\leq j\leq r}})\not\equiv0$.

We define the  Cartier divisor $(\Db,\o)\subsetneq(Z,\o)$ by the principal ideal
\beq\label{Eq.Discrim.Via.Classical.Discrim}
(\Db,\o):=\{ \frD(\{a^{(j)}_{m_1,\dots,m_k}\}_{\substack{\suml_i m_i=p_j\\1\leq j\leq r}})=0\}\sset (Z,\o).
\eeq

By the definition of $\frD$:  $\o'\in\Db$ if and only if $(X,\o')$  is
not ordinary along $(Z,\o')$.
Note that $\frD$ is a polynomial, therefore this construction ``preserves the category'': if the germs $(Z,\o)\sset(X,\o)$ are algebraic/analytic/formal/etc.
 then so is the subgerm $(\De^\bot,\o)\sset(Z,\o)$.

\

This definition can be restated more geometrically as follows. The choice of a good basis of $I_{(X,\o)}$ near $(Z,\o)$, defines a
 rational map from $(Z,\o)$ to the parameter space of complete intersections:
\beq\label{Eq.Map.of.Sing.Locus.to.Parameter.Space.of.Complete.Intersections}
(Z,\o)\ni\o'\stackrel{\phi}{\to}\Big(\{a^{(1)}_{m_1,\dots,m_k}(\o')\}_{\substack{\suml_i m_i=p_1}},\dots,
\{a^{(r)}_{m_1,\dots,m_k}(\o')\}_{\substack{\suml_i m_i=p_r}}\Big)\in\prodl_{j=1}^r|\cO_{\P^{k-1}}(p_j)|.
\eeq
As $(X,\o)$ is generically ordinary along $(Z,\o)$ this map is generically well defined. Its indeterminacy locus consists of those points $\o'\in Z$ where at
least one of the collections of coefficients vanishes, i.e. the multiplicity of some $f_i$ jumps. The discriminant of transversal type is
 the pullback: $\Db=\phi^*(\De)$.

\bprop\label{Thm.Discriminant.via.Classical.Discrim}
The defining ideal $I_{(\Db,\o)}\sset\cO_{(Z,\o)}$ is independent
 of all the choices made (the local coordinates in $(\k^{n+r},\o)$, the basis of $I_{(Z,\o)}$, the good basis of $I_{(X,\o)}$).
\eprop
Indeed, as is proved in proposition \ref{Thm.Normal.Cone.is.well.defined}, the change of bases/coordinates results
in the action $\P GL_{(\cO_{(Z,\o)})}(k)\circlearrowright(Z,\o)\times \P^{k-1}$. This action is linear, it preserves the classical discriminant.
 Thus it does not change the defining ideal of $(\Db,\o)\sset(Z,\o)$.

\bex\label{Ex.Discriminant.Eq.for.generic.transversal.A1}
In the simplest case let $(X,\o)=\{f(\ux)=0\}\sset(\k^{n+1},\o)$ be a hypersurface singularity,  with $(Z,\o)=Sing(X,\o)$ a complete intersection.
 Suppose the generic transversal type of $(X,\o)$ along $(Z,\o)$ is ordinary of multiplicity two. (In Siersma's notations this is called: $A_1$-transversal type.)
Then $f=\suml^k_{i,j=1}a_{ij}g_ig_j$, and we can assume that the matrix $\{a_{ij}\}_{ij}$ is symmetric. The discriminant
is then  $\Db=\{det\Big(\{a_{ij}\}_{ij}|_{(Z,\o)}\Big)=0\}$.
 Suppose $(Z,\o)$ is smooth and $dim(Z,\o)=1$, i.e. $k=n$, then the generic singularity type of $(X,\o)$ along $(Z,\o)$ is $A_\infty$.
 Take a deformation of $(X,\o)$ that preserves $Sing(X,\o)$ and splits $\De^\bot$ into a few reduced points.
 Near such points the local equation of $(X,\o)$ can be brought to the form $\{\suml^{n-2}_{i=1}x^2_i+x_{n-1}x^2_n=0\}$,
   the standard notation for this singularity type is $D_\infty$. Then we get: the number of these $D_\infty$ points is the degree of the scheme
 $\{det\Big(\{a_{ij}\}_{ij}|_{(Z,\o)}\Big)=0\}$. This recovers  \cite[theorem 7.18]{Pellikaan-PhD}, see also \cite[page 176]{de Jong}.
\eex
\bex\label{Ex.Singular.Locus.is.non.smooth.Transversal.Type.Ordinary} Let $(X,\o)=\{(xy)^pa=z^pb\}\sset(\k^3,0)$,
where $a,b\in \cO_{(\k^3,0)}$ are invertible.
 Here $I_{(Z,\o)}=(xy,z)$ and $\P(\cN_{(X,\o)/(Z,\o)})=\{y^p_1a|_{(Z,\o)}=y_2^pb|_{(Z,\o)}\}\sset\P^1\times(Z,\o)$
  is a hypersurface smooth over $(Z,\o)$. Thus $\De^\bot$ is not supported at the origin, eventhough $(Z,\o)$ is singular. Note also that the flat deformation
$\{(xy+t)^pa=z^pb\}\sset(\k^3,0)\times(\k^1,\o)$ induces a (flat) smoothing $\cZ=\{Z_t\}$ of $(Z,\o)$, while preserving the generic vanishing order.
And for $t\neq0$  all the fibres of $\P\cN_{X_t/Z_t}$ are smooth, thus $\De^\bot_{t\neq0}=\empty$. Compare to the flatness of $\Db$ in deformations,
 proposition \ref{Thm.Discriminant.Flat.Deformations}.
\eex

\subsubsection{The global case}
Suppose we begin from quasi-projective or analytic (for $\k=\C$) spaces $Z\sset X$. Then the local/pointwise definition of germs $(\De^\bot,\o)\sset(Z,\o)$
globalizes. Here the germs are algebraic/analytic (i.e. all the local rings are either localizations of affine or analytic), thus we can take
 the representatives/open neighborhoods and glue along them.
\bprop
The local divisors $\{(\Db,\o)\sset(Z,\o)\}_{\o\in Z}$ glue to the global effective Cartier divisor $\De^\bot\sset Z$.
\eprop
\bpr
The defining ideal of each germ $(\Db,\o)$ is  principal. Thus it is enough to prove that these ideals glue to a coherent sheaf of ideals,
 $I_\Db\sset\cO_Z$. Namely, we should check compatibility: given a germ $(Z,\o)$ with some representatives $\cU_1\into\cU_2$,
  the identification of sheaves $\cO_{\cU_1}\isom i^*\cO_{\cU_2}$ induces the identity isomorphism $I_{\De^\bot}(\cU_1)\isom i^*(I_{\De^\bot}(\cU_2))$.
  And this follows as $\De^\bot$ does not depend on the choice of coordinates/representatives/bases of ideals,
   see Proposition \ref{Thm.Discriminant.via.Classical.Discrim}.
\epr

\subsection{The defining ideal of $(\Db,\o)\sset (Z,\o)$}\label{Sec.Db.Definition.via.Fitting.Ideals}
The discriminant of transversal type is defined in the last section as the pullback of the classical discriminant. It is often useful to work directly with the
 ideal $I_{(\De^\bot,\o)}\sset\cO_{(Z,\o)}$ or the sheaf $I_{\De^\bot}\sset\cO_{Z}$. These are directly obtained using
  \S\ref{Sec.Discriminant.Classical.Critical.Locus.Explicit}.

Fix the complete intersections $Sing(X,\o)=(Z,\o)\sset(X,\o)$, suppose $(X,\o)$ is s.c.i. over $(Z,\o)$.
Fix a basis  $I_{(Z,\o)}=(g_1,\dots,g_k)$ and a good basis $I_{(X,\o)}=(f_1,\dots,f_r)$.
 Then $\P\cN_{(X,\o)/(Z,\o)}=\{\tf_1=\cdots=\tf_r=0\}\sset(Z,\o)\times\P^{k-1}$,
  where $\tf_j\in \cO_{(Z,\o)}[y_1,\dots,y_k]$ is obtained as the leading term of $f_j$. The projection $\P\cN_{(Z,\o)}\to(Z,\o)$ is precisely the
   projection $\pi$ of \S\ref{Sec.Discriminant.Classical}. Thus equation \eqref{Eq.Critical.Locus.Global}
   and definition \ref{Def.Classical.Discriminant} give us:
\bcor\label{Thm.Discriminant.Via.Fitting.Ideals}
1. The critical locus of the projection is the subscheme:
\[Crit(\pi):=\Big\{\tf_1=\cdots=\tf_r=0,\ \ rank\Big( d \tf_1,\dots,d \tf_r\Big)<r\Big\}\sset\P^{k-1}\times(Z,\o),\]
where  $\{d\tf_i\}_i$ are $k\times 1$ columns of partial derivatives of $\{\tf_i\}_i$, taken with respect to the homogeneous coordinates in $\P^{k-1}$.
\\2. Suppose the restriction $Crit(\pi)\stackrel{\pi}{\to}(Z,\o)$ is a finite map. Then the defining ideal of $\Db\sset(Z,\o)$ is: $Fitt_0(\pi_*\cO_{Crit(\pi)})$.
\ecor

\subsection{Further stratifications of the singular locus}\label{Sec.Db.Stratification.Singular.Locus}
 Recall that at some points of $(Z,\o)=Sing(X,\o)$ the ($\mu=const$) singularity type
 of $(L^\bot\cap X,\o)$ depends on the choice of the section $L^\bot$, see example
 \ref{Ex.Transversal.Type.Depends.choice.of.Section}.
Therefore we make the stratification according to the singularities of the fibers of the projectivized normal cone,
  $\P\cN_{(X,\o)/(Z,\o)}$.

Any stratification of the parameter space, $|\cO_{\P^{k-1}}(p)|$ or $\prodl_i|\cO_{\P^{k-1}}(p_i)|$, e.g. by singularity type  for
 some equivalence relation, induces a stratification of $Sing(X)$. More precisely, using the map $\phi$ of equation
 \eqref{Eq.Map.of.Sing.Locus.to.Parameter.Space.of.Complete.Intersections}, we get the following:

\begin{center} if (the closure of) some stratum   $\overline\Sigma\sset\prodl_i|\cO_{\P^{k-1}}(p_i)|$   is defined by an ideal $I_{\overline\Sigma}$ then

the ideal $\phi^*(I_{\overline\Sigma})\sset\cO_Z$ defines the corresponding stratum on the singular locus.
\end{center}
\bex Consider the $\mu=const$  stratification of $|\cO_{\P^{k-1}}(p)|$: the points of a stratum correspond to all the hypersurfaces
 that can be deformed to a given hypersurface in a $\mu=const$ way.  This  defines the strata:
 \beq
Sing(X)\supset\De^\bot=\lSi_{A_1}\supset\lSi_{A_2},\lSi_{A_1,A_1}\supset\lSi_{A_3},\lSi_{A_1,A_2},\lSi_{A_1A_1A_1}\supset \dots\lSi_{D_4}\supset\cdots
\eeq
\eex

\section{Some general properties of $\Db$}\label{Sec.Db.General.Properties}
The definition of $\Db$ as the pullback of the classical discriminant is somewhat theoretical, as in most cases it is
extremely difficult to write down the classical discriminant
explicitly. (Recall that even in the hypersurface case, $r=1$,
$\frD$ is a polynomial of degree $k(p-1)^{k-1}$ in $\bin{p+k}{k}$
variables.) Also the computation of the Fitting ideal $Fitt_0(\pi_*\cO_{Crit(\pi)})$ is, in general, an involved procedure.
Yet, some consequences are obtained immediately.

\subsection{The discriminant pulls back}\label{Sec.Db.Gen.Prop.Pulls.Back}
Suppose we are given morphisms of (germs of) manifolds, as on the diagram.
\parbox{13cm}{
Here $X_1=\phi^*(X_2)$ and $Z_1=\phi^*(Z_2)$ are pullback of schemes/ideals.
 Assume $X_i$, $Z_i$ are reduced, l.c.i., $Z_i$ is a connected component of $Sing(X_i)$ and $X_i$ is s.c.i. over $Z_i$.
 Suppose, moreover, $X_i$ is generically ordinary along $Z_i$ and the multiplicity sequences, of $X_1$ along $Z_1$ and $X_2$ along $Z_2$, coincide.
\bprop
   Then $\Db_{X_1/Z_1}=\phi^*\Db_{X_2/Z_2}$.
\eprop
\bpr It is enough to check the statement locally at each point. Thus we work with germs. We have:
 $I_{(Z_1,\o_1)}=\phi^*(I_{(Z_2,\o_2)})=\bl \phi^*(g_1),\dots,\phi^*(g_k)\br$ and $I_{(X_1,\o_1)}=\phi^*(I_{(X_2,\o_2)})=\bl \phi^*(f_1),\dots,\phi^*(f_r)\br$,
  and in both cases the sequences are regular.
}\quad\quad
$\bM M_1&\stackrel{\phi}{\to}&M_2\\\cup&&\cup\\X_1=\phi^*(X_2)&\stackrel{\phi}{\to}& X_2\\\cup&&\cup\\Z_1=\phi^*(Z_2)&\stackrel{\phi}{\to}& Z_2
\\\cup&&\cup\\\Db_1&&\Db_2\eM$

 To define $\Db_{(X_1,\o_1)/(Z_1,\o_1)}$ we expand:
\beq
\phi^*(f_i)=\suml_{\sum m_j=p_i}\phi^*(g^{m_1}_1)\cdots \phi^*(g^{m_k}_k) \tilde{a}^{(i)}_{m_1\dots m_k}.
\eeq
 But we can also pullback the initial expansions,
 $f_i=\suml_{\sum m_j=p_i}g^{m_1}_1\cdots g^{m_k}_k a^{(i)}_{m_1\dots m_k}$.
 As $(\phi^*(g_1),\dots,\phi^*(g_k))$ form a regular sequence, we get $\tilde{a}^{(i)}_{m_1\dots m_k}=\phi^*a^{(i)}_{m_1\dots m_k}$.
 Therefore:
\begin{multline}
I_{\Db_{(X_1,\o_1)/(Z_1,\o_1)}}=\Big(\frD(\{\tilde{a}^{(i)}_{m_1,\dots,m_k}\}_{\substack{\suml_i m_i=p_i\\1\leq i\leq r}})\Big)=
\Big(\frD(\{\phi^*a^{(i)}_{m_1,\dots,m_k}\}_{\substack{\suml_i m_i=p_i\\1\leq i\leq r}})\Big)=\\=
 \Big(\phi^*\frD(\{a^{(i)}_{m_1,\dots,m_k}\}_{\substack{\suml_i m_i=p_i\\1\leq i\leq r}})\Big)=
\phi^*I_{\Db_{(X_2,\o_2)/(Z_2,\o_2)}}.  \quad\bull
\end{multline}
%\text{\epr}
\bex
Consider the surface $X_1=\{x^2z^q=y^2+x^3\}\sset M_1\approx\k^3$, cf. example \ref{Ex.Whitney.Umbrella}. Then $X_1$ is the
 pullback of $X_2=\{x^2z=y^2+x^3\}\sset M_2\approx\k^3$, under the covering $(x,y,z)\stackrel{\phi}{\to}(x,y,z^q)$. Thus
\beq
\Db_{X_1/Z_1}=\phi^*(\Db_{X_2/Z_2})=\{z^q=0\}\sset Z_1=\{x=0=y\}.
\eeq
\eex

\beR\label{Ex.Conditions.on.Pulback}
(Importance of $\phi^*(Z)$ being reduced.)
 Let  $M'=\k^n\stackrel{\phi}{\to} M=\k^n$ by $(y_1,\dots,y_n)\to(y^2_1,\dots,y_n)$. Then for $X=\{x^2_1=x^2_2\}$ we have $Z=\{x_1=0=x_2\}$ and
 $\Db_{X/Z}=\varnothing$. But $\phi^*(Z)=\{y^2_1=0=y_2\}$ and $\phi^*(X)=\{y^4_1=y^2_2\}$,
 the transversal type is generically non-ordinary, i.e. $\Db_{X'/Z'}=Z'$.
\eeR

\subsection{The discriminant is determined by infinitesimal neighborhood of $Z$ in $X$}\label{Sec.Db.Gen.Prop.Defined.by.Infinit.Neighb}
By its construction $\De^\bot$ reflects degenerations of the projectivized normal cone and does not depend on the degenerations of
\\\parbox{14.6cm}{
'higher order terms'.
 This ideas is made precise by a variation of the last proposition. Fix two triples (of germs)  $\Big\{Z_i=Sing(X_i)\sset X_i\sset(M_i)\Big\}_{i=1,2}$.
  Suppose $X_i$  is generically ordinary along $Z_i$ and the multiplicity sequences in both cases are the same: $p_1\le p_2\le\cdots\le p_r$.
\bprop\label{Thm.Db.is.Determined.by.infinit.Neighb}
Suppose the restriction $Z_1\stackrel{\phi|_{Z_1}}{\to}Z_2$ is an isomorphism and moreover
 $\phi^*(I_{X_2}\otimes  \quotients{\cO_{M_2}}{I^{p_r+1}_{Z_2}})=
I_{X_1}\otimes \quotients{\cO_{M_1}}{I^{p_r+1}_{Z_1}}$. Then $\phi^*(\De^\bot_2)=(\De^\bot_1)$.
\eprop
 }\quad
  $\bM M_1&\stackrel{\phi}{\to}& M_2\\\cup&\underset{\phi|_{Z_1}}{\ \ }&\cup\\ Z_1 &\isom & Z_2\\\cup&&\cup\\ \De^\bot_1 && \De^\bot_2\eM$

\bpr As before, it is enough to check the statement pointwise.
As $(\De^\bot,\o)$ is fully determined by $\P\cN_{(X,\o)/(Z,\o)}$, it is enough to show that  $\phi$ induces
 isomorphism $\P\cN_{\quotients{(X_1,\o)}{(Z_1,\o)}}\isom\P\cN_{\quotients{(X_2,\o)}{(Z_2,\o)}}$.
 Fix a basis $I_{(Z_2,\o)}=(g_1,\dots,g_k)$ so that $I_{(Z_1,\o)}=(\phi^*(g_1),\dots,\phi^*(g_k))$. Fix a good basis
  $I_{(X_2,\o)}=(f_1,\dots,f_r)$ then the assumption reads: there exist $\tau_1,\dots,\tau_r\in I_{(Z_1,\o)}^{p_r+1}$ such that
   $(\phi^*(f_1)+\tau_1,\dots,\phi^*(f_r)+\tau_r)$ is a basis of $I_{(X_1,\o)}$. It follows that this is a good basis. But then the
   expansion $f_i=\suml_{\sum m_j=p}g^{m_1}_1\dots g^{m_k}_k a^{(i)}_{m_1\dots m_k}\in \cO_{(M_2,\o)}$ ensures the expansion:
\beq
\phi^*(f_i)=\suml_{\sum m_j=p_i}\phi^*(g^{m_1}_1)\dots \phi^*(g^{m_k}_k)\big(\phi^*(a^{(i)}_{m_1\dots m_k})+b^{(i)}_{m_1\dots m_k}\big)\in\cO_{(M_1,\o)},
 \text{ for some }
  b^{(i)}_{m_1\dots m_k}\in I_{(Z_1,\o)}.
  \eeq
 Thus
 \beq
 \P\cN_{(X_1,\o)/(Z_1,\o)}=\{\suml_{\sum m_j=p_i}y^{m_1}_1\dots y^{m_k}_k \phi^*(a^{(i)}_{m_1\dots m_k})=0,\ \forall\ i\}\isom \P\cN_{(X_2,\o)/(Z_2,\o)}.\text{\epr}
 \eeq
The proposition states that $\Db$ is determined by the $(p_r+1)$-infinitesimal neighborhood of $(Z,\o)$ in $(X,\o)$.
 Therefore $\Db$ is determined by the formal neighborhood:
\bcor
Given two triples $(M_1,X_1,Z_1)$ and $(M_2,X_2,Z_2)$ with $Z_i=Sing(X_i)$.
 Suppose  $Z_1\approx Z_2$ and the completions along the singular loci are isomorphic
$\widehat{(M_1,X_1)}\approx\widehat{(M_2,X_2)}$. Then the discriminants are (embedded) isomorphic.
\ecor

The converse statement to proposition \ref{Thm.Db.is.Determined.by.infinit.Neighb} does not hold: if the map $M_1\stackrel{\phi}{\to}M_2$
 restricts to an isomorphism $\phi|_Z$, with $\phi^*(\Db_2)=\Db_1$, this  does not imply much relation of $(X_1,\o)$ to $(X_2,\o)$.
 For example, compare $I_{(X_1,\o)}=(x(zx^3+y^3)+x^5+y^5)$ and  $I_{(X_2,\o)}=(z(x^4+y^4)+x^2y^2+x^5+y^5)$.
 In both cases $(Z,\o)=Sing(X_i,\o)=\{x=0=y\}$
  and their generic type along $(Z,\o)$ is the ordinary multiple point of multiplicity 4.
   For $(X_1,\o)$ the degeneration of transversal type at $\o$ is: 3 roots collide to a triple root. For $(X_2,\o)$ the degeneration is:
    two pairs of roots collide to two double roots. Thus in both cases $\De^\bot=(z^2=0)$.

\subsection{Flat deformations}\label{Sec.Db.Gen.Prop.Deforms.Flatly}
We prove that $\De^\bot$ deforms flatly in those flat deformations of $X$ that induce flat deformations of (reduced!) singular locus $Z=Sing(X)$
 and preserve the multiplicity sequence. More precisely, given a good basis $I_{(X,\o)}=(f_1(\ux),\dots,f_r(\ux))$, fix a flat deformation
  $I_{(\cX,\o)}=(f_1(\ux,t),\dots,f_r(\ux,t))$ with $f_i(\ux,0)=f_i(\ux)$, such that the (reduced) singular locus $\cZ=Sing(\cX)$ is a flat family:
   $I_{(\cZ,\o)}=(g_1(\ux,t),..,g_k(\ux,t))$ and   $I_{(Z,\o)}=(g_1(\ux,0),..,g_k(\ux,0))$.

\bprop\label{Thm.Discriminant.Flat.Deformations}
 Suppose $(X,\o)$ is s.c.i. over $(Z,\o)$
and $(\cX,\o)$ is s.c.i. over $(\cZ,\o)$ and the multiplicity sequence is preserved.
 Then the family $\De^\bot_{\cX/\cZ}$ is flat and its central fibre is $\De^\bot_{X/Z}$.
\eprop
\bpr
By the assumption we can use the standard expansion
$f_j(\ux,t)=\sum g_\uI(\ux,t)^{m_\uI}a_{m_\uI}(\ux,t)$. Thus $\Db_{\cX/\cZ}=\{\frD(a_{\uI}(\ux,t))=0\}$
 is a flat family  that  specializes to $\De^\bot_{X/Z}$. (Note that $\frD(a_{\uI}(\ux,t))$ is a
 power series in $t$.)
\epr
In many cases this property allows the quick computation of the transversal multiplicity of $\Db$.

\

Example \ref{Ex.Whitney.Umbrella} shows that $\De^\bot$ can be non-reduced if the degeneration occurs
`faster than normally'. Another reason for being non-reduced is when the degeneration is not `minimal'.
\bex\label{Ex.Whitney.Umbrelly.Fast.Dying.Mult=p}
Consider the surface $X=\{x^pz=y^p+x^{p+1}\}\sset\k^3$. Its singular locus is the line $Z=\{x=0=y\}$.
Consider the projection $\k^3\proj Z$, $(x,y,z)\to z$, and the fibres $\pi^{-1}(z)$.
Then we have a family of plane curve singularities, $\pi^{-1}(t)\cap X\sset\pi^{-1}(t)=(\k^2,\o)$, for $t\in Z$.
 This family is equimultiple, thus the projectivized
tangent cones of these curve singularities form the flat family: $\{\si^p_xz=\si^p_y\}\sset\P^1_{\si_x,\si_y}\times\k^1_z$.
For each $z\neq0$
there are $p$ distinct roots, while for $z=0$ all these roots coincide, thus $\Db$ is supported at $z=0$.
Now the multiplicity can be computed using a flat deformation or via the critical locus and the Fitting ideal.
 \bee[i.]
\item
Under a generic deformation
this multiple root at $z=0$ splits into several double roots near $z=0$. In our case  one can take
$\{\si^p_xz=\si^p_y-\ep\si^{p-1}_x\si_y\}\sset\P^1_{\si_x,\si_y}\times\k^1_z$. By direct check, for each fixed $\ep\neq0$
the number of the double roots near $z=0$ is $(p-1)$. So, by the flatness of $\Db$ in deformations,  the multiplicity of $\Db$ for the initial surface is $(p-1)$.
\item
Blow-up $\k^3$ along the line $Z=\{x=0=y\}$, let $E\sset Bl_Z\k^3$ be the exceptional divisor, consider the strict
transform $\tX\sset Bl_Z\k^3$ and the projection $\tX\cap E\to Z$. Explicitly:
$Z\times\P^1=E\supset \tX\cap E=\{\si^p_xz=\si^p_y\}\stackrel{\pi}{\to} \{z\}\sset Z$. This is
a $p:1$ covering, totally ramified over $z=0$ and the ramification degree is $(p-1)$.
 The critical locus is (see proposition \ref{Thm.Discriminant.Via.Fitting.Ideals}) $\{\si^{p-1}_xz=0=\si^{p-1}_y\}$.
  Thus $\cO_{Crit(\pi)}\approx\quotients{\k[[z,\si_y]]}{(z,\si^{p-1}_y)}$ and $\pi_*(\cO_{Crit(\pi)})\approx\big(\quotients{\k[z]}{(z)}\big)^{p-1}$,
 as a module over $\cO_{(Z,\o)}=\k[[z]]$. Therefore $Fitt_0(\pi_*\cO_{Crit(\pi)})=(z^{p-1})$.
 \eee
\eex
\bex\label{Ex.Discrim.Depends.on.The.Tangent.Cone.Only}
Consider the hypersurface $X=\{z(x^p+y^p)=x^{q}+y^{r}\}\sset\k^3$, with $p<q,r$. Again, $Sing(X)=\{x=0=y\}$
and the discriminant is a point on $z$-axis, namely, $\Db=\{x=y=z=0\}$, as a set. The deformation
 $\{zx^p+(z-\ep)y^p=x^{q}+y^{r}\}$ splits the discriminant
into two: at $z=0$ and at $z=\ep$. The previous example gives that both points have
multiplicity $(p-1)$, {\em regardless} of $q,r$. (Compare to proposition \ref{Thm.Db.is.Determined.by.infinit.Neighb}.) Hence the multiplicity in the current case is $2(p-1)$.
\eex

If in the family $\{X_t\}$ the generic multiplicity along $(Z_t,\o)$ changes then the (non-flat) family $\Db_t$ is not semicontinuous in any reasonable sense. For example, consider  $(X_t,\o)=\{x^p+y^p+t(x^2+zy^2)\}\sset(\k^3,\o)$. Then
 $\Db_{t=0}=\varnothing$, while  $\Db_{t\neq0}=\{(0,0,0)\}$.

\subsection{Comparison of $\Db$ to L\^{e} numbers/cycles}\label{Sec.Db.vs.Le.Numbers}
We give just one example to show that the relation to L\^{e} numbers of Massey is not at all obvious. We work in the notations of \cite[Chapter 1]{Massey-book}.
\bex
For $(X,\o)=\{y^p=\frac{x^{p+1}}{p+1}+\frac{t^qx^p}{p}\}$ the ideal defining the singular scheme is $(y^{p-1},t^{q-1}x^p,x^p+t^qx^{p-1})$. Its saturation gives $\Si_f=V(y^{p-1},x^{p-1})$.

Then the polar schemes are $\Ga^2=V(y^{p-1})$ and $\Ga^1=V(y^{p-1},x^p,x^p+t^qx^{p-1})\smin\Si_f=V(y^{p-1},x+t^q)$.

Now the L\^{e} cycles are:
\beq\ber
\La^1=[(y^{p-1},x^p,x^p+t^qx^{p-1})]-[(y^{p-1},x+t^q)]=[(y^{p-1},x^{p-1})]\quad and\\ \La^0=[(y^{p-1},x+t^q,t^{q-1}x^p)]=(q-1+pq)(p-1)[(y,x,t)].
\eer\eeq
 Thus $\la^1=dim_\k\quotients{\k[x,y,t]}{(t,x^{p-1},y^{p-1})}=(p-1)^2$, while  $\la^0=(q-1+pq)(p-1)$.

But $deg[\Db]=q(p-1)$, as can be seen e.g. by deformation $(X_t,\o)=\{y^p=\frac{x^{p+1}}{p+1}+\frac{(t^q-\ep)x^p}{p}\}$, see also example \ref{Ex.Whitney.Umbrelly.Fast.Dying.Mult=p.Recomputed}.
\eex

\subsection{The transversal multiplicity of the discriminant}\label{Sec.Db.Multiplicity}
Proposition \ref{Thm.Classical.Discriminant.Multiplicity} implies the following formula: if  $\pi^{-1}(\o)$ has only isolated singularities and $(Z,\o)$
 contains the miniversal deformation of the singularities of $\pi^{-1}(\o)$ then
\beq
mult(\De^\bot,\o)=\mu_{total}(\pi^{-1}(\o))+\mu'_{total}(\pi^{-1}(\o)).
\eeq
  In most cases of interest
  $(Z,\o)=Sing(X,\o)$ is of low dimension and does not contain the miniversal deformation of $\pi^{-1}(\o)$. Thus the only conclusion is
   $mult(\De^\bot,\o)\ge\mu_{total}(\pi^{-1}(\o))+\mu'_{total}(\pi^{-1}(\o))$, see example \ref{Ex.Whitney.Umbrella}.
 (Indeed, we have the usual map from $(Z,\o)$ to the germ of miniversal deformation, such that the family is the pullback. And
  computation of $mult(\Db,\o)$ corresponds to the intersection of the discriminant in the miniversal deformation by a not-necessarily-transversal,
  not-necessarily-smooth curve-germ.)

\bprop
In the notations of \S\ref{Sec.Db.Definition.via.Fitting.Ideals} suppose the map $Crit(\pi)\stackrel{\pi}{\to}(Z,\o)$ is finite at $\o\in Z$.

1. Let $\{pt_\al\}\sset Crit(\pi)$ be the points of the fibre $\pi^{-1}(\o)$, let $\De_\al^\bot\sset(Z,\o)$ be the corresponding discriminants.
Then, as Cartier divisors, $\Db=\suml_\al \De_\al^\bot$.

2. Take a  one dimensional complete intersection subgerm $(C,\o)\subseteq(Z,\o)$,
 such that  $(C,\o)\cap\Db\sset(C,\o)$ is zero dimensional (Cartier divisor).
Then $deg\Big((C,\o)\cap\Db\Big)=deg\Big(Crit(\pi)\cap\pi^{-1}(C,\o)\Big)$.

3. Suppose $(X,\o)\sset(\k^{n+1},\o)$ is a hypersurface and $(Z,\o)$ is smooth. Suppose near each $pt_\al\in Crit(\pi)$ the defining equation of $\tX\cap E$
 in some local coordinates has the form $f_\al(\ux)+h_\al(\uz)=0$, here $\uz$ are the local coordinates of $(Z,\o)$.
  Then $mult(\Db,\o)=\suml_{pt_\al\in Crit(\pi)}\mu(f_\al(\ux))\cdot ord(h_\al(\uz))$.
\eprop
In part 2 on both sides we have the scheme-theoretic intersection, defined by the union of the ideals.
By the degree of a zero-dimensional scheme we mean the length of its ring: $deg(Y)=dim_\k\cO_Y$.
\\\bpr
{\bf 1.} In this case $Crit(\pi)$ is a multi-germ, thus $\pi_*(\cO_{Crit(\pi)})=\oplusl_\al \pi_*(\cO_{(Crit(\pi),pt_\al)})$ .
 For the direct sum of modules one has: $Fitt_0(M_1\oplus M_2)=Fitt_0(M_1)\cdot Fitt_0(M_2)$.
 Thus $I_{\Db_{}}=\prodl_\al I_{\Db_\al}$.

{\bf 2.} We should prove the two equalities: $deg\Big((C,\o)\cap\Db\Big)=deg(\De^\bot_{\pi|_C})=deg\Big(Crit(\pi)\cap\pi^{-1}(C,\o)\Big)$.

The left equality is immediate by base-change, as $\Db$ is defined by pulling back the classical discriminant, \S\ref{Sec.Db.Definition.via.Classical}.

Therefore we can restrict to $(C,\o)$. So,
 we  assume that $(Z,\o)$ is a one-dimensional locally complete intersection and $\Db\sset(Z,\o)$
 is a Cartier divisor (in particular it is a zero dimensional subscheme).
 We should prove: $deg(\Db)=deg(Crit(\pi))$.

 By part 1 it is enough to consider the one-point critical locus, $pt=Crit(\pi)\in \tX\cap E$.
\bei
\item We start from the case: $(Z,\o)$ is smooth. Note that
\begin{multline}
deg(\Db)=dim_\k\cO_{\Db}=dim_\k\quotient{\cO_{(Z,\o)}}{Fitt_0(\pi_*\cO_{Crit(\pi)})}\quad
\\ and\quad
deg(Crit(\pi))=dim_\k\cO_{Crit(\pi)}=dim_\k(\pi_*\cO_{Crit(\pi)}).
\end{multline}
Thus, the statement to prove is: given a finite module $M$ over a one-dimensional
regular local ring, $\cO_{(Z,\o)}$, the colength of the Fitting ideal satisfies: $colength_{\cO_{(Z,\o)}}(Fitt_0M)=dim_\k M$. This is a standard statement of commutative algebra. Take the minimal
free resolution: $\cO_{(Z,\o)}^{\oplus p}\stackrel{A}{\to}\cO_{(Z,\o)}^{\oplus q}\to M\to0$. As $M$ is finite,
it is supported at one point only, so $p\ge q$. Furthermore, as the ring is local and regular, $A$ is equivalent,
by $A\to UAV$, to a diagonal matrix. Let $z$
be a generator of $\cO_{(Z,\o)}$, then $Fitt_0(M)=Fitt_q(A)=(z^{\sum d_i})$, here $\{d_i\}$ are the
exponents of the diagonal. Thus $M\approx\oplus_i \cO_{(Z,\o)}/(z^{d_i})$ and $colength(Fitt_0M)=\sum d_i=dim_\k M$.
 Proving that $deg(\Db)=deg(Crit(\pi))$.

\item Suppose $(Z,\o)$ is a complete intersection (of dimension one), then it can be smoothed.
Let $\{Z_t\}_{t\in(\k^1,\o)}$
be a smoothing, then we have the (flat) family of projections, $(\tX_t\cap E)\proj Z_t$.
 Explicitly, if $\tX\cap E=\{\{f_\al(\ux,\uz)=0\}_\al\}$ then $\tX_t\cap E=\{\{f_\al(\ux,\uz_t)=0\}_\al\}$.

This induces
the flat family $\{\Db_t\sset Z_t\}$. Thus, for $t\in(\k^1,\o)$ small enough, we can fix some
(small enough, Zariski open) neighborhood of $\o\in Z$ such that $deg(\Db_{t=0})=deg(\Db_{t})$. Here the r.h.s.
is the total degree of $\Db_t$ in the neighborhood. Note that $\Db_{t\neq0}$ is a subscheme of a smooth curve $Z_t$.
 Thus the statement holds for $\Db_{t\neq0}$ and then, by flatness, for $\Db_{t=0}$.
 \eei

{\bf 3.} Restrict from $(Z,\o)$ to the generic curve $(C,\o)\sset(Z,\o)$, then
\beq
mult(\Db,\o)=length(\pi_*\cO_{Crit(\pi)})=length(\cO_{Crit(\pi)})=
\suml_{pt_\al\in Crit(\pi)}length(\cO_{(Crit(\pi),pt_i)}).
\eeq
Thus it is enough to verify the claim for each point. In the local coordinates:
\beq
(Crit(\pi),pt_\al)=\{\di_{\ux} f(\ux)=0=f(\ux)+h(z)\}.
\eeq
Working locally, we can redefine $z$ to get $h(z)=z^p$, where $p=ord(h)$. Then,  as a $\k[\ux]$-module,
$\cO_{(Crit(\pi),pt_i)}\approx \quotients{\k[\ux]}{\di_{\ux}f}\cdot<1,z,\dots,z^{p-2}>$. Therefore $length(\cO_{(Crit(\pi),pt_\al)})=\mu(f_\al)\cdot p$.
\epr
\bex
\bee[i.]
\item Using part 3 of the  proposition one gets immediately the multiplicity
 of the discriminant in examples
 \ref{Ex.Whitney.Umbrella},  \ref{Ex.Whitney.Umbrelly.Fast.Dying.Mult=p},
 \ref{Ex.Discrim.Depends.on.The.Tangent.Cone.Only}.
 For example \ref{Ex.Singular.Locus.is.non.smooth.Transversal.Type.Ordinary} (with non-smooth $(Z,\o)$) one can use Part 2.

\item  \label{Ex.Whitney.Umbrelly.Fast.Dying.Mult=p.Recomputed}
 (Extending example \ref{Ex.Whitney.Umbrelly.Fast.Dying.Mult=p}.)
 Consider the hypersurface singularity $X=\{x^p_1x^q_n=f_p(x_2,\dots,x_{n-1})+g_{>p}(x_1,\dots,x_n)\}\sset(\k^n,\o)$,
where $f_p(x_2,\dots,x_{n-1})$ is a homogeneous form of degree $p$, while  $g_{>p}(x_1,\dots,x_n\in\cm^{p+1})$. Suppose $f_p$ is generic,
so that $\{f_p(x_2,\dots,x_{n-1})=0\}\sset\P^{n-3}$ is smooth. Suppose $g_{>p}$ contains
a monomial $x^N_1$ for some $N$. Then $Sing(X)=\{x_1=\cdots=x_{n-1}=0\}\sset\k^n$ and the
generic transversal type (for $x_n\neq0$) is ordinary. The discriminant $\Db\sset Sing(X)$ is supported
at the point $\{x_n=0\}\sset Sing(X)$ and its multiplicity equals the length
of the scheme
$Crit(\pi)=\{\si^{p-1}_1x^q_n=0=\dv f_p(\si_2,\dots,\si_{n-1})\}\sset\k^1_{x_n}\times\P^{n-2}_{\si_1,\dots,\si_{n-1}}$.
As the form $f_p$ is generic, this scheme coincides with the scheme $\{x^q_n=0=\dv f_p(\si_2,\dots,\si_{n-1})\}$,
whose degree is $q(p-1)^{n-2}$.

\item Consider the hypersurface singularity
$(X,\o)=\{x^p+x^{p-r}y^rg_1(\uz)+y^pg_2(\uz)+h_{>p}(x,y,\uz)=0\}\sset(\k^2_{x,y},\o)\times(Z,\o)$.
Here $g_{>p}(x,y,\uz)\in (x,y)^{p+1}$.
The (reduced) singular locus is of codimension two: $(Z,\o)=\{x=0=y\}$. The strict transform under the blows-up along $(Z,\o)$
 is: $\tX\cap E=\{\si_x^p+\si_x^{p-r}\si_y^rg_1+\si_y^pg_2=0\}\sset\P^1_{\si_x,\si_y}\times(Z,\o)$.
The critical set is:
$Crit(\pi)=\{p\si^{p-1}_x+(p-r)\si^{p-r-1}_x\si^r_y g_1=0=r\si^{p-r}_x\si^{r-1}_yg_1+p\si^{p-1}_yg_2=\si_y^pg_2\}$.
 Therefore $\si^p_x\in I_{Crit(\pi)}$, hence $Crit(\pi)$ is located in the $\si_y=1$ chart of $\P^1$.
By part 2 of the proposition, $mult(\Db,\o)=length(Crit(\pi)|_{(C,\o)})$, for any generic curve germ $(C,\o)\sset(Z,\o)$.
 Let $ord(g_1(\uz))=q_1$ and $ord(g_2(\uz))=q_2$, then
\bei
\item if $q_1\ge q_2$ then $mult(\Db,\o)=q_2(p-1)$.
\item if $q_1< q_2$ then $mult(\Db,\o)=q_1(p-1)+(p-r)(q_2-q_1)$.
\eei
\eee
\eex

\bex
Take two hypersurface  germs, $(X_i,\o)=\{f_i=0\}\sset(\k^{n+1},\o)$. Suppose they share the singular locus, $Sing(X_1,\o)=(Z,\o)=Sing(X_2,\o)$,
  and $(Z,o)$ is
 a reduced complete intersection of dimension $n-1$.
 %Suppose both $(X_i,\o)$ are s.c.i. over $(Z,\o)$ and
 %are generically ordinary along $(Z,\o)$.
  In this case the projectivized normal cones, $\P\cN_{(X_1,\o)/(Z,\o)}$, $\P\cN_{(X_2,\o)/(Z,\o)}$, are subschemes of $\P^1\times(Z,\o)$.
  Suppose these subschemes are disjoint. Then $\Db_{f_1f_2}=\Db_{f_1}+\Db_{f_2}$, i.e. $I_{(\Db_{f_1f_2},\o)}=I_{(\Db_{f_1},\o)}I_{(\Db_{f_2},\o)}$.
\eex

\subsection{The $\cO_{(Z,\o)}$-resolution of $\cO_{(\Db,\o)}$ and defining equation of $(\Db,\o)$}\label{Sec.Db.Resolution.Defining.Equation.Tangent.Cone}
We want to obtain some information about the local equation of the discriminant. By the previous proposition, if $Crit(\pi)\cap\pi^{-1}(\o)$
contains several points, then $\Db$ is the sum of the components. Thus it is enough to consider the case of one critical point, i.e.
the map of  germs $(Crit(\pi),pt)\to(\Db,\o)$.
 Let $(\cX,pt)=\{\underline\tf(\ux)+\ua(\ux,\uz)=0\}\sset (\P\cN_{\quotients{(X,\o)}{(Z,\o)}},pt)$, where $\tf_i\in\k[[\ux]]$ while $a_i\in \cm_{(Z,\o)}[[\ux]]$, i.e. $\ua(\ux,0)=0$, and $\ua(0,\uz)\neq0$.
Thus, we can consider $\{\underline\tf=0\}$ as an isolated singularity and $\ua$ as its deformation.

Note that we work locally in $(\P\cN_{\quotients{(X,\o)}{(Z,\o)}},pt)$, thus the singularity is not generically ordinary in any sense.

\subsubsection{A free resolution of $\pi_*\cO_{(Crit(\pi),pt)}$ as an $\cO_{(Z,o)}$-module.}
\bthe
Suppose both $(\tf_1,\dots,\tf_r)$ and $(\tf_2,\dots,\tf_r)$ define isolated (complete intersection) singularities in $(\k^{n+r},\o)$.
Fix their Milnor numbers, $\mu=\mu(\tf_1,\dots,\tf_r)$, $\mu_{\hat{1}}=\mu(\tf_2,\dots,\tf_r)$, as in \S \ref{Sec.Discriminant.Classical.Tranversal.Multiplicity}.
 Then the free resolution is
 \[0\to \cO^{\oplus(\mu+\mu_{\hat{1}})}_{(Z,\o)}\stackrel{[\tf_1+a_1]}{\to}\cO^{\oplus(\mu+\mu_{\hat{1}})}_{(Z,\o)}\to \pi_*(\cO_{(Crit(\pi),pt)})\to0.\]
(The linear map $[\tf_1+a_1]$ is defined during the proof.)

 Thus the locally defining ideal $I_{(\Db,\o)}\sset\cO_{(Z,\o)}$ of  $\Db$ is generated by $det[\tf_1+a_1]$.
\ethe
\bpr
{\bf The hypersurface case} (In this cases we follow  e.g. \cite{Teissier1976}.)

{\bf Step 1.} By corollary \ref{Thm.Discriminant.Via.Fitting.Ideals} the critical locus can be presented in the form:
\beq
\cO_{(Crit(\pi),pt)}=\quotients{\cO_{(Z,\o)}[[\ux]]}{(\tf+a,\di_1(\tf+a),\dots,\di_k(\tf+a))},
\eeq
the derivatives are taken with respect to $\ux$ coordinates.

Choose some $\k$-basis of the the Milnor algebra $\quotients{\k[[\ux]]}{(\di_1\tf,\dots,\di_k\tf)}$, and fix some $\k[[\ux]]$-representatives of this basis,
 $\{v_\al\}$.
 Use the composition of maps $\k[[\ux]]\into \cO_{(Z,\o)}[[\ux]]\to \cO_{(Crit(\pi),pt)}$ and denote the images of $v_\al$ by $[v_\al]$.

We claim that $\{[v_\al]\}_\al$ generate $\pi_*(\cO_{(Crit(\pi),pt)})$, i.e. generate $\cO_{(Crit(\pi),pt)}$ as an $\cO_{(Z,\o)}$ module.
 Indeed, we have
\beq
 \pi_*(\cO_{(Crit(\pi),pt)})\sseteq Span_{\cO_{(Z,\o)}}\{[v_\al]\}+\cm_{(Z,\o)}\pi_*(\cO_{(Crit(\pi),pt)}).
 \eeq
Take the quotient of $\pi_*(\cO_{(Crit(\pi),pt)})$ by $Span_{\cO_{(Z,\o)}}\{v_\al\}$ and apply the Nakayama lemma.

\

{\bf Step 2.}
It remains to understand the relations among $\{[v_\al]\}$, i.e. the kernel of the surjection
\beq
\cO^{\oplus\mu}_{(Z,\o)}\to \pi_*(\cO_{(Crit(\pi),pt)})\to0.
\eeq
The relations come from the ideal $(\tf+a,\di_1(\tf+a),\dots,\di_k(\tf+a))$.
 We should pass from the $\cO_{(Z,\o)}[[\ux]]$-module structure on $\cO_{(Crit(\pi),pt)}$ to the $\cO_{(Z,\o)}$-module structure.
  In other words, we should express the action of $\ux$  via $\cO_{(Z,\o)}$-action. By the construction of $\{v_\al\}$:
\beq
x_iv_\al=\suml_\be b^{(i,\al)}_{\be} v_\be+\suml_j c^{(i,\al)}_{j} \di_j \tf\in\k[[\ux]],\quad \text{ where } b^{(i,\al)}_{\be}\in\k,\quad  c^{(i,\al)}_{j}\in \k[[\ux]].
\eeq
Here the coefficients $\{b^{(i,\al)}_{\be}\}$ are defined uniquely, while the coefficients $\{c^{(i,\al)}_{j}\}$ are unique up to
 the Koszul relations of the regular sequence $\di_1\tf,\dots,\di_k\tf$.
 Therefore in $\cO_{(Crit(\pi),pt)}$ we have:
$[x_i][v_\al]=\suml_\be b^{(i,\al)}_{\be}[ v_\be]+\suml_j [-c^{(i,\al)}_{j}\di_j a]$. In general  $c^{(i,\al)}_{j}\di_j a$ can still depend on $\ux$,
 more precisely - they lie in $\cm_{(Z,\o)}[[\ux]]$. Expand them again, in terms of $\{v_\al\}$ and $\{\di_j f\}$, and iterate the procedure, until one gets:
\beq
 x_i[v_\al]=\suml_\be b^{(i,\al)}_{\be}(\uz)[ v_\be]\in\cO_{Crit(\pi)},\quad \text{where now } b^{(i,\al)}_{\be}(\uz)\in\cO_{(Z,\o)}.
\eeq
Despite the (Koszulian) non-uniqueness of $\{c^{(i,\al)}_{j}\}$ the resulting coefficients $\{b^{(i,\al)}_{\be}(\uz)\}\in\cO_{(Z,\o)}$ are well defined, as any
 Koszul relation among $\{\di_j \tf\}$ induces that among $\{\di_j(\tf+a)\}$.

\

{\bf Step 3.}
The last equation defines the needed action $x_i\circlearrowright \cO_{(Crit(\pi),pt)}$, we denote the corresponding $\cO_{(Z,\o)}$-linear operator by $[x_i]$.
Accordingly, for any element $h\in\cO_{(Z,\o)}[[\ux]]$ we have the operator $[h]\circlearrowright \cO_{(Crit(\pi),pt)}$.
 Thus we have the operators$[\tf+a]$, $[\di_1(\tf+a)]$, \dots, $[\di_k(\tf+a)]$. By the construction: $[\di_i(\tf+a)]=0$ for all $i$.
 (Here it is important to notice that the ideal $(\di_1(\tf),\dots,\di_k(\tf))$ is a complete intersection,
  the only relations among its generators are Koszul and thus extend to the relations among $(\di_1(\tf+a),\dots,\di_k(\tf+a))$.)
Thus the   only possibly non-trivial operator is $[\tf+a]=0$. Thus we get the presentation:
\beq
\cO^{\oplus\mu}_{(Z,\o)}\stackrel{[\tf+a]}{\to}\cO^{\oplus \mu}_{(Z,\o)}\to \pi_*(\cO_{(Crit(\pi),pt)})\to0.
\eeq

Finally we note that the module $\pi_*\cO_{Crit(\pi)}$  is a torsion, of rank=0, because $\ua(0,z)\neq0$.
 Therefore the map $[\tf+a]$ cannot have a kernel, thus is injective, i.e. we have
\beq
0\to \cO^{\oplus\mu}_{(Z,\o)}\stackrel{[\tf+a]}{\to}\cO^{\oplus \mu}_{(Z,\o)}\to \pi_*(\cO_{(Crit(\pi),pt)})\to0.
\eeq

\

{\bf The case of complete intersections.}

By corollary \ref{Thm.Discriminant.Via.Fitting.Ideals} the critical locus can be presented in the form:
\beq
\cO_{(Crit(\pi),pt)}=\quotient{\cO_{(Z,\o)}[[\ux]]}{\Big(\{\tf_i+a_i\}_{i=2,\dots, r},\ Fitt_0(\di_j(\tf_i+a_i))\Big)},
\eeq
  where the derivatives are taken with respect to $\ux$ coordinates.

As in the hypersurface case, choose some $\k$-basis of the vector space $\quotients{\k[[\ux]]}{(\tf_2,\dots,\tf_r,\ Fitt_0\bpm \di_j \tf_i\epm)}$,
 and fix some $\k[[\ux]]$-representatives of this basis,
 $\{v_\al\}$.  The size of the basis is $\mu(\tf_1,\dots,\tf_r)+\mu(\tf_2,\dots,\tf_r)$, by L\^{e}-Greuel formula
   (see \S\ref{Sec.Discriminant.Classical.Tranversal.Multiplicity}).
 Use the composition $\k[[\ux]]\into \cO_{(Z,\o)}[[\ux]]\to \cO_{(Crit(\pi),pt)}$ and denote the images of $v_\al$ by $[v_\al]$.

\

As in the hypersurface case one gets: $\{[v_\al]\}$ generate $\pi_*(\cO_{(Crit(\pi),pt)})$,  as an $\cO_{(Z,\o)}$-module.

\

It remains to understand the relations, i.e. the kernel of the surjection
\beq
\cO^{\oplus(\mu+\mu_{\hat{1}})}_{(Z,\o)}\to \pi_*(\cO_{(Crit(\pi),pt)})\to0.
\eeq
The relations come from the ideal $\Big(\{\tf_i+a_i\},\ Fitt_0(\di_j(\tf_i+a_i))\Big)$.

We begin with the part  $\Big(\tf_2+a_2,\dots,\tf_r+a_r,\ Fitt_0(\di_j(\tf_i+a_i))\Big)$.
Unlike the hypersurface case, this ideal is not a complete intersection.
 Yet, the scheme $V\Big(\tf_2+a_2,\dots,\tf_r+a_r,\ Fitt_0(\di_j(\tf_i+a_i))\Big)\sset (Z,\o)\times Spec(\k[\ux])$ is a flat family of schemes over $(Z,\o)$.
  Therefore any syzygy of $\Big(\tf_2,\dots,\tf_r,\ Fitt_0(\di_j \tf_i)\Big)$ extends to that of  $\Big(\tf_2+a_2,\dots,\tf_r+a_r,\ Fitt_0(\di_j(\tf_i+a_i))\Big)$.

 As in the hypersurface case,
  we should express the action of $\ux$  via $\cO_{(Z,\o)}$-action, i.e. should define the operators $[x_i]\circlearrowright \cO_{(Crit(\pi),pt)}$.

   By the construction of $\{v_\al\}$, as in the hypersurface case:
\beq
x_i v_\al=\suml_\be b^{(i,\al)}_\be v_\be+\suml_\ga c^{(i,\al)}_\ga h_\ga\in\k[[\ux]],\quad  \text{ where } b^{(i,\al)}_\be\in\k, \quad  c^{(i,\al)}_\ga\in\k[[\ux]],
\eeq
  while $\{h_\ga\}$ are the generators of $(\tf_2,\dots,\tf_r,Fitt_0(\di_j \tf_i)$. Here the coefficients $\{b^{(i,\al)}_\be\}$ are defined uniquely,
   while $\{c^{(i,\al)}_\ga\}$ are defined up to the (not necessarily Koszul) relations in $(\tf_2,\dots,\tf_r,Fitt_0(\di_j \tf_i)$.
Therefore in $\cO_{(Crit(\pi),pt)}$ we have:
$x_i[v_\al]=\suml_\be b^{(i,\al)}_\be [v_\be]+\suml_\ga [c^{(i,\al)}_\ga h_\ga]$. As in the hypersurface case $[c^{i,\al}_\ga]$ might still depend on $\ux$,
 more precisely: $[c^{i,\al}_\ga]\in\cm_{(Z,\o)}[[\ux]]$. Iterate the procedure until one  gets:
\beq
x_i [v_\al]=\suml_\be b^{(i,\al)}_\be(\uz) [v_\be],\quad  \text{ where } b^{(i,\al)}_\be(\uz)\in\cO_{(Z,\o)}.
\eeq
Despite the intermediate non-uniqueness, due to relations among $\{h_\ga\}$, the resulting coefficients $\{b^{(i,\al)}_\be\}\in\cO_{(Z,\o)}$ are unique.
 (Again, because of the flatness, any relation among $\tf_2,\dots,\tf_r,\ Fitt_0(\di_j(\tf_i))$ extends to a relation among
  $\tf_2+a_2,\dots,\tf_r+a_r,\ Fitt_0(\di_j(\tf_i+a_i))$.)

\

Thus we have the needed (well defined) action $x_i\circlearrowright \cO_{(Crit(\pi),pt)}$, we denote the corresponding linear operator by $[x_i]$.
Accordingly, for any element $h\in\cO_{(Z,\o)}[[\ux]]$ we have the operator $[h]\circlearrowright \cO_{(Crit(\pi),pt)}$.

As in the hypersurface case, we get the tautology: $[\tf_2+a_2]=0$,\dots, $[\tf_r+a_r]=0$,  $[Fitt_0(\di_j (\tf_i+a_i))]=0$.
   (Here the flatness of $V\Big(\tf_2+a_2,\dots,\tf_r+a_r,\ Fitt_0(\di_j(\tf_i+a_i))\Big)$ is used.)
The only relation to understand is  $[\tf_1+g_1]=0$, this operator is in general non-trivial.
 Thus we get the presentation of the $\cO_{(Z,\o)}$-module $\pi_*(\cO_{(Crit(\pi),pt)})$:
\beq
\cO^{\oplus(\mu+\mu_{\hat{1}})}_{(Z,\o)}\stackrel{[\tf_1+a_1]}{\to}\cO^{\oplus(\mu+\mu_{\hat{1}})}_{(Z,\o)}\to \pi_*(\cO_{(Crit(\pi),pt)})\to0
\eeq

Finally, as in the hypersurface case, we note that the module $\pi_*\cO_{Crit(\pi)}$  is a torsion, of rank=0, because $\ua(0,z)\neq0$.
 Therefore the map $[\tf_1+a_1]$ cannot have a kernel, thus is injective.
\epr

\beR
The proposed resolution is certainly not the only possible. In fact, in the hypersurface case, if $\tf$ is not weighted homogeneous, i.e.
 $[\tf]\neq0 \in \quotients{\k[[\ux]]}{(\di_i \tf)}$,  then $\tf$ can be taken as one of $\{v_\al\}$. Then the matrix $[\tf+a]$ contains an entry
  invertible in $\cO_{(Z,\o)}$,  i.e. the resolution is non-minimal.

   In the hypersurface case one could start from a basis of the
   Tjurina algebra, $\quotients{\k[[\ux]]}{(\tf,\di_1 \tf,\dots,\di_k \tf)}$. Lift it to $\k[[\ux]]$ and send to $\cO_{(Crit(\pi),pt)}$. Then, as before, we get:
    $\cO_{(Z,\o)}^{\oplus \tau}\to \pi_*\cO_{(Crit(\pi),pt)}\to0$. However, the relations are more complicated now, they come from syzygies of the Tjurina
     algebra (which is not a complete intersection).
  And, unlike the hypersurface case, the definition of $[x_i]\circlearrowright\pi_*\cO_{(Crit(\pi),pt)}$ is more complicated now, as the family
 $V(\tf+a,\di_1(\tf+a),\dots,\di_k(\tf+a))\sset (Z,\o)\times Spec(\k[[\ux]])$ is not flat.
\eeR

 Though the defining equation, $(\Db,\o)=\{det[\tf_1+a_1]=0\}\sset(Z,\o)$,   cannot be written in an explicit form, below
 we can get some information on the participating monomials.

\subsubsection{The discriminant for deformations by constant terms}
Suppose in some local coordinates the universal family is $\cX=\{\tf_1(\ux)+s_1=0,\dots,\tf_r(\ux)+s_r=0\}\sset(\k^{k},\o)\times Spec(\k[\{s_i\}])$.
 Here $\{s_i\}$ are independent variables, the case of functions on the singular locus, $\{a_i(\uz)\}\in\cO_{(Z,\o)}$, is obtained by base change.
 (Recall, that the singularity $\{\tf_i\}$ is not generically ordinary.)

Denote by $\mu$ the Milnor number of the ICIS $\{\tf_1=\cdots=\tf_r=0\}\sset(\k^{k},\o)$.
 Suppose for any $1\le j\le r$ the ideal $I_\hj:=(\tf_1,\dots,\tf_{j-1},\tf_{j+1},\dots,\tf_r)$ defines an isolated (complete intersection) singularity.
  Define the auxiliary Milnor number, $\mu_{\hj}:=\mu(I_\hj)$, as in \S\ref{Sec.Discriminant.Classical.Tranversal.Multiplicity}.
\bcor
Then   the defining equation of $\Db\sset (Z,\o)$ has the form: $\suml^r_{j=1} c_j s^{\mu+\mu_\hj}_j+(\{s_is_j\}_{j\neq i})=0$.
\ecor
Here $\{c_i\in \k^*\}$ are some non-zero constants, while $(\{s_is_j\}_{j\neq i})$ is a collection of monomials involving at least two distinct $s_i$'s.
\bpr By the assumption, $det[\tf_1+s_1]$ is a polynomial in $s_1,\dots,s_r$.
 Note that  $\big(det[\tf_1+s_1]\big)|_{s_i=0}=det\big([\tf_1+s_1]|_{s_i=0}\big)$.
  Put $s_2=\cdots=s_r=0$, then, by remark \ref{Thm.Classical.Discriminant.Multiplicity.one.dimensional.base},
   $det\big([\tf_1+s_1]|_{s_2=\cdots=s_r=0}\big)=s^{\mu+\mu_{\hat{1}}}_1$. Similarly,
  \beq
  det([\tf_1+s_1]|_{\{s_i=0\}_{i\neq j}})=s^{\mu+\mu_\hi}_i.
\eeq
Hence the statement.
\epr

\subsubsection{The discriminant for the weighted homogeneous case}
Consider the family
\beq
\cX=\{\tf_1(\ux)+h_1(\ux,\us^{(1)})=0,\dots,\tf_r(\ux)+h_r(\ux,\us^{(r)})=0\}\sset(\k^k,\o)\times Spec(\k[\{\us^{(j)}\}]).
\eeq
Here, using multi-indices, $h_j(\ux,\us^{(j)})=\suml_{\um} s^{(j)}_\um\ux^\um$.

\bprop\label{Thm.Total.Degree.Disriminant}
 Suppose $\tf_1,\dots,\tf_r\in\k[[\ux]]$ are weighted homogeneous, of degree $\{w(\tf_j)\}$ with respect to the weights $\{w(x_i)\}$.
 Then $\Db$ is weighted
homogeneous, and the only possible monomials
 to appear are:
 $\prodl_j\prodl_{\um} (s^{(j)}_{\um})^{n^{(j)}_{\um}}$, where
 \[
 \suml^r_{j=1}\suml_\um n^{(j)}_{\um} \Big(w(\tf_j)-\suml_i m_i w(x_i)\Big)=
 (\mu+\mu_{\hat{1}})w(\tf_1)=\suml_{j=1}^rw(\tf_j)\prodl^k_{i=1}\Big(\frac{w(\tf_j)}{w(x_i)}-1\Big)\prodl^r_{\substack{q=1\\q\neq j}}\frac{w(\tf_q)}{w(\tf_j)-w(\tf_q)}.
 \]
\eprop
\bpr
Impose the condition ``$(\tf_j+h_j)$ is weighted homogeneous, of weight $w(\tf_j)$", then
 the weights of the coefficients are  fixed, $w(s^{(j)}_{m_I})=w(\tf_j)-\suml_i m_i w(x_i)$.

  Then $\Db$ is weighted homogeneous. We know that the monomials $\{(s^{(j)}_{0\dots 0})^{\mu+\mu_\hj}\}$ are present.
  Thus the total (weighted) degree of $\Db$ is $(\mu+\mu_\hj)\cdot w(\tf_j)$. Thus the only possible monomials in the defining equation of $\Db$ are:
\beq
\prodl_j\prodl_{\um} (\us^{(j)}_{\um})^{n^{(j)}_{\um}},\quad \text{ where }\quad
 \suml ^r_{j=1}\suml_\um n^{(j)}_{\um} \Big(w(\tf_j)-\suml_i m_iw(x_i)\Big)=
 (\mu+\mu_{\hat{1}})w(\tf_1).
\eeq
It remains to compute $(\mu+\mu_{\hat{1}})$.
We use the expression for Poincar\'{e} series of weighted homogeneous complete intersection with isolated singularity, of dimension $n$ and codimension $r$
 \cite[Satz 3.1]{Greuel-Hamm}:
 \beq
P_{\tf_1\dots \tf_r}(t)=res_{\tau=0}\frac{\tau^{-n-1}}{1+\tau}\Big[\prodl^{n+r}_{i=1}\frac{1+\tau t^{w(x_i)}}{1-t^{w(x_i)}}
\prodl^r_{j=1}\frac{1-t^{w(\tf_j)}}{1+\tau t^{w(\tf_j)}}
+1\Big],\quad\quad
\mu=P_{\tf_1\dots \tf_r}(1).
 \eeq
So, we should extract the residue and take the limit $t\to1$.

Consider here $\{w(\tf_j)\}$ as independent variables in $\R_{>0}$, then $P_{\tf_1\dots \tf_r}(t)$ depends continuously on $\{w(\tf_j)\}$.
 Therefore we can compute under the assumptions: $\{w(\tf_j)\neq w(\tf_i)\}_{i\neq j}$ and $\{w(\tf_j)\neq0\}$. After the computation is done, the cases
 $w(\tf_j)=w(\tf_i)$  are obtained by taking the limit.

The expression $R(t,\tau):=\frac{\tau^{-n-1}}{1+\tau}\Big[\prodl^{n+r}_{i=1}\frac{1+\tau t^{w(x_i)}}{1-t^{w(x_i)}}
\prodl^r_{j=1}\frac{1-t^{w(\tf_j)}}{1+\tau t^{w(\tf_j)}}\Big]$ is a rational function in $\tau$, with poles at $\tau=0,-1,\{-t^{-w(f_j)}\}$.
 For $\tau\to\infty$ this function decreases as $\frac{1}{\tau^2}$. Therefore
 \beq
 res_{\tau=0}(R(t,\tau))+ res_{\tau=-1}(R(t,\tau))+\suml^r_{j=1} res_{\tau=-t^{-w(f_j)}}(R(t,\tau))=0.
 \eeq
Assuming $w(\tf_j)\neq0$, we have $res_{\tau=-1}(R(t,\tau))=(-1)^{n}2$. Assuming $w(\tf_j)\neq w(\tf_i)$ we get:
\beq
res_{\tau=-t^{-w(\tf_j)}}(R(t,\tau))=
\frac{(-t^{-w(\tf_j)})^{-n-1}}{1-t^{-w(\tf_j)}}  \prodl^k_{i=1}\frac{1-t^{w(x_i)-w(\tf_j)}}{1-t^{w(x_i)}}
\Big(\prodl^r_{\substack{q=1\\q\neq j}}
\frac{1-t^{w(\tf_q)}}{1-t^{w(\tf_q)-w(\tf_j)}}\Big)\frac{1-t^{w(\tf_j)}}{t^{w(\tf_j)}}.
\eeq
Therefore
\beq
\liml_{t\to1}\Big(res_{\tau=-t^{-w(\tf_j)}}(R(t,\tau))\Big)=
-\prodl^k_{i=1}\Big(\frac{w(\tf_j)}{w(x_i)}-1\Big)
\Big(\prodl^r_{\substack{q=1\\q\neq j}}\frac{w(\tf_q)}{w(\tf_j)-w(\tf_q)}\Big).
\eeq
Altogether we get:
\beq
\mu=P_{\tf_1\dots \tf_r}(1)=(-1)^{n-1}2+\suml_{j=1}^r\prodl^k_{i=1}\Big(\frac{w(\tf_j)}{w(x_i)}-1\Big)
\Big(\prodl^r_{\substack{q=1\\q\neq j}}\frac{w(\tf_q)}{w(\tf_j)-w(\tf_q)}\Big).
\eeq
Similarly:
\beq
\mu_{\hat{1}}=P_{\tf_2\dots \tf_r}(1)=(-1)^{n}2+\suml_{j=2}^r\prodl^k_{i=1}\Big(\frac{w(\tf_j)}{w(x_i)}-1\Big)
\Big(\prodl^r_{\substack{q=2\\q\neq j}}\frac{w(f_q)}{w(\tf_j)-w(\tf_q)}\Big).
\eeq
Finally, combine these two to get:
\begin{multline}\label{Eq.Proof.last.equation}
\mu+\mu_{\hat{1}}=
\suml_{j=2}^r\prodl^k_{i=1}\Big(\frac{w(\tf_j)}{w(x_i)}-1\Big)
\Big(\prodl^r_{\substack{q=2\\q\neq j}}\frac{w(\tf_q)}{w(\tf_j)-w(\tf_q)}\Big)\frac{w(\tf_j)}{w(\tf_j)-w(\tf_1)}
+
\prodl^k_{i=1}\Big(\frac{w(\tf_1)}{w(x_i)}-1\Big)\Big(\prodl^r_{\substack{q=2}}\frac{w(\tf_q)}{w(\tf_1)-w(\tf_q)}\Big)=
\\=
\suml_{j=1}^r\frac{w(\tf_j)}{w(\tf_1)}\prodl^k_{i=1}\Big(\frac{w(\tf_j)}{w(x_i)}-1\Big)\prodl^r_{\substack{q=1\\q\neq j}}\frac{w(\tf_q)}{w(\tf_j)-w(\tf_q)}.
\end{multline}
This finishes the proof.
\epr

\vspace{-0.5cm}
\bex
In the hypersurface case, $r=1$, equation \eqref{Eq.Proof.last.equation} gives the  Milnor number of a weighted homogeneous isolated hypersurface singularity,
  $\mu=\prodl_i\big(\frac{w(\tf)}{w(x_i)}-1\big)$, cf. \cite{Milnor-Orlik}. Thus the necessary condition for a monomial
   $\prodl_{\um} (\us_{\um})^{n_{\um}}$ to participate in $\Db$ is:
\beq
\suml_\um n_{\um}\big(1-\suml_i m_i\frac{ w(x_i)}{w(\tf)}\big)=\prodl_i\big(\frac{w(\tf)}{w(x_i)}-1\big).
\eeq
For example, let $\tf=\sum x^{p_i}_i$  and $h(\ux,\uz)=\suml_{i=1}^k s_i x_i+s_0$. Then the only possible monomials are:
 $s^{n_0}_0s_1^{n_1}\cdots s^{n_r}_r$, where
 $n_0+\sum n_i(1-\frac{1}{p_i})=\prod(p_i-1)$.
 \eex


\begin{thebibliography}{99}
\bibitem[Aluffi-1995]{Aluffi-1995} P. Aluffi, {\em Singular schemes of hypersurfaces},
Duke Math. J. 80 (1995), no. 2, 325--351

%\bibitem[Aluffi-1996]{Aluffi-1996} P. Aluffi, {\em Chern classes for singular hypersurfaces},



 \bibitem[Aluffi-2005]{Aluffi-2005} P. Aluffi {\em Characteristic classes of singular varieties.}
  Topics in cohomological studies of algebraic varieties, 1--32, Trends Math., Birkh\"{a}user, Basel, 2005.


%\bibitem[Aluffi-Faber]{Aluffi-Faber} P. Aluffi, E. Faber {\em Splayed divisors and their Chern classes}

\bibitem[AGLV-book1]{AGLV} V.I. Arnol'd, V.V. Goryunov, O.V. Lyashko, V.A. Vasil'ev, {\em
Singularity theory.I.} Reprint of the original English edition
from the series Encyclopaedia of Mathematical Sciences [Dynamical
systems. VI, Encyclopaedia Math. Sci., 6, Springer, Berlin, 1993].
Springer-Verlag, Berlin, 1998.

\bibitem[AGLV-book2]{AGLV2} V.I. Arnol'd, V.A. Vasil'ev, V.V. Goryunov, O.V. Lyashko, {\em Singularities. II.
 Classification and applications.} (Russian) With the collaboration of B. Z. Shapiro.
 Itogi Nauki i Tekhniki, Current problems in mathematics. Fundamental directions, Vol. 39, 5256, Akad. Nauk SSSR,
 1989.


\bibitem[Bennett-1977]{Benett-1977} B. Bennett, {\em Normally flat deformations.}
Trans. Amer. Math. Soc. 225 (1977), 1--57


\bibitem[Benoist2012]{Benoist} O. Benoist, {\em Degr\'{e}s d'homog\'{e}n\'{e}it\'{e} de l'ensemble des intersections
  compl\`{e}tes singuli\`{e}res.}  Ann. Inst. Fourier (Grenoble) 62 (2012), no. 3, 1189--1214.


\bibitem[C.C.D.R.S.2011]{C.C.D.R.S.2011} E. Cattani, M.A. Cueto, A. Dickenstein, S.Di Rocco, B.  Sturmfels,
 {\em Mixed Discriminants}, arXiv:1112.1012


\bibitem[Dimca-book]{Dimca92} A. Dimca, {\it Singularities and topology of hypersurfaces.} Universitext. Springer-Verlag,
New York, 1992.

\bibitem[Eisenbud-book]{Eisenbud-book} D.  Eisenbud, {\em Commutative algebra. With a view toward algebraic geometry.}
 Graduate Texts in Mathematics, 150. Springer-Verlag, New York, 1995.

\bibitem[Esterov2011]{Esterov} A. Esterov {\em The discriminant of a system of equations}. Adv. Math. 245 (2013), 534--572.



\bibitem[Fulton]{Fulton-book} W.  Fulton, {\em Intersection theory}. Second edition.
A Series of Modern Surveys in Mathematics, 2. Springer-Verlag, Berlin, 1998.

\bibitem[G.K.Z.]{Gelfand-Kapranov-Zelevinsky} I.M. Gelfand, M.M. Kapranov, A.V. Zelevinsky,
{\em Discriminants, resultants and multidimensional determinants.}
 Reprint of the 1994 edition. Modern Birkh\"{a}user Classics. Birkh\"{a}user Boston, Inc., Boston, MA, 2008. x+523 pp.

\bibitem[Greuel]{Greuel} G. M. Greuel, {\em Der Gauss-Manin Zusammenhang isolierter Singularit\"{a}äten von vollst\"{a}ändigen Durchschnitten},
  Dissertation, G\"{o}öttingen, 1973.

\bibitem[Greuel-Hamm]{Greuel-Hamm} G.-M. Greuel, H.A. Hamm, {\em Invarianten quasihomogener vollst\"{a}ändiger Durchschnitte.}
Invent. Math. 49 (1978), no. 1, 67--86.


\bibitem[Heinzer-Kim-Ulrich]{Heinzer-Kim-Ulrich}
W. Heinzer, M.-K. Kim, B. Ulrich, {\em The Gorenstein and complete intersection properties of associated graded rings.}
 J. Pure Appl. Algebra 201 (2005), no. 1--3, 264--283

\bibitem[Heitmann-Jorgensen]{Heitmann-Jorgensen} R.C. Heitmann, D.A.Jorgensen {\em Are complete intersections
complete intersections?} J. Algebra 371 (2012), 276--299.


\bibitem[de Jong1990]{de Jong} Th.de Jong, {\em The virtual number of $D_\infty$ points I.}
 Topology 29 (1990), no. 2, 175--184
\bibitem[de Jong-de Jong1990]{de Jong-de Jong} J.de Jong, Th.de Jong, {\em The virtual number of $D_\infty$ points II.}
 Topology 29 (1990), no. 2, 185--188

%\bibitem[Kazarian2000]{Kazarian2000} M. Kazarian, {\em Classifying spaces of singularities and Thom polynomials.}
% New developments in singularity theory (Cambridge, 2000), 117--134,
%NATO Sci. Ser. II Math. Phys. Chem., 21, Kluwer Acad. Publ., Dordrecht, 2001

%\bibitem[Kazarian2003]{Kazarian2006} M. Kazarian, {\em Thom polynomials.} Singularity theory and its applications,
%85--135, Adv. Stud. Pure Math., 43, Math. Soc. Japan, Tokyo, 2006

\bibitem[K.K.N.]{Kazarian.Kerner.Nemethi} M. Kazarian, D. Kerner, A. N\'{e}methi, {\em
    Discriminant of the ordinary transversal singularity type. The global properties},  arXiv:1308.6045.

%\bibitem[Kerner-Nemethi-2012]{Durfee2} D. Kerner, A. Nemethi, {\em The 'corrected' Durfee's inequality for
%homogeneous complete intersections}, to appear in Mathematische Zeitschrift.

\bibitem[L\^{e}]{Le} D. T. L\^{e}, {\em Calculation of Milnor number of isolated singularity of complete intersection}
, Funct. Anal. Appl. 8 (1974), 127--131.


\bibitem[Looijenga-book]{Looijenga-book} E. Looijenga {\em Isolated Singular Points on Complete Intersections}
 London Math. Soc. LNS 77, CUP, 1984

\bibitem[Massey]{Massey-book} D.B. Massey, {\em L\^{e} cycles and hypersurface singularities.} Lecture Notes in Mathematics, 1615. Springer-Verlag, Berlin, 1995. xii+131 pp

\bibitem[Milnor-Orlik]{Milnor-Orlik} J. Milnor, P. Orlik, {\em  Isolated singularities defined by weighted homogeneous polynomials.}
Topology 9 (1970), 385--393.

\bibitem[Mond-Pellikaan]{Mond-Pellikaan} D. Mond, R. Pellikaan, {\em
Fitting ideals and multiple points of analytic mappings.} Algebraic geometry and complex analysis (Pátzcuaro, 1987), 107–161,
Lecture Notes in Math., 1414, Springer, Berlin, 1989

\bibitem[N\u{a}st\u{a}sescu-Van Oystaeyen]{Nastasescu-Van Oystaeyen}
C. N\u{a}st\u{a}sescu, F. Van Oystaeyen, {\em Graded and filtered rings and modules.}
Lecture Notes in Mathematics, 758. Springer, Berlin, 1979. x+148 pp. ISBN: 3-540-09708-2


\bibitem[Pellikaan1985]{Pellikaan-PhD} G.R. Pellikaan, {\em Hypersurface singularities and resolutions of Jacobi modules.}
 Dissertation, Rijksuniversiteit te Utrecht, Utrecht, 1985. Drukkerij Elinkwijk B. V., Utrecht, 1985. vii+168 pp

\bibitem[Pellikaan1990]{Pellikaan} R.Pellikaan {\em Deformations of hypersurfaces with a one-dimensional singular locus.}
J. Pure Appl. Algebra 67 (1990), no. 1, 49--71


 \bibitem[Piene1977]{Piene1978} R. Piene, {\em Some formulas for a surface in $\P^3$.} Algebraic geometry
 (Proc. Sympos., Univ. Troms{\o}, Troms{\o}, 1977), pp. 196--235,
 Lecture Notes in Math., 687, Springer, Berlin, 1978

\bibitem[Piene1978]{Piene1977} R. Piene, {\em Polar classes of singular varieties.} Ann. Sci. \'{E}cole Norm. Sup.
(4) 11 (1978), no. 2, 247--276.

\bibitem[Seade-book]{Seade} J. Seade, {\it On the Topology of Isolated Singularities in Analytic Spaces.} Progress in
Mathematics 241, Birkh\"auser 2006.

\bibitem[Siersma1987]{Siersma1987} D. Siersma, {\em Singularities with critical locus a 1-dimensional complete intersection and transversal type $A_1$.}
 Topology Appl. 27 (1987), no. 1, 51--73.

\bibitem[Siersma2000]{Siersma2000} D. Siersma, {\em The vanishing topology of non isolated singularities.} New developments
in singularity theory (Cambridge, 2000), 447--472, NATO Sci. Ser. II Math. Phys. Chem., 21, Kluwer Acad. Publ.,
Dordrecht, 2001

\bibitem[van Straten2011]{van Straten} D.van Straten {\em Gorenstein-duality for one-dimensional almost complete
intersections-with an application to non-isolated real singularities}, arXiv:1104.3070

\bibitem[Teissier1976]{Teissier1976} B. Teissier, {\em The hunting of invariants in the geometry of discriminants.}
 Real and complex singularities (Proc. Ninth Nordic Summer School/NAVF Sympos. Math., Oslo, 1976), pp. 565--678.
 Sijthoff and Noordhoff, Alphen aan den Rijn, 1977.


\bibitem[Valabrega-Valla]{Valabrega-Valla} P. Valabrega, G. Valla, {\em Form rings and regular sequences.}
Nagoya Math. J. 72 (1978), 93--101.

\bibitem[Vasconcelos]{Vasconcelos} W.-V. Vasconcelos, {\em Computational methods in commutative algebra and algebraic geometry.}
  With chapters by David Eisenbud, Daniel R. Grayson, J\"{u}ürgen Herzog and Michael Stillman. Algorithms and Computation in Mathematics, 2.
   Springer-Verlag, Berlin, 1998. xii+394 pp. ISBN: 3-540-60520-7
\end{thebibliography}
\end{document}